\documentclass[12pt,a4paper]{amsart}
\usepackage[english]{babel}
\usepackage{amssymb,latexsym,amsfonts,amsthm,upref,amsmath}
\usepackage[margin=1in]{geometry}
\usepackage[colorlinks=true]{hyperref}
\hypersetup{citecolor=blue, urlcolor=blue, colorlinks=false, linkcolor=red, citecolor=green, filecolor=magenta}

\newtheoremstyle{prim}{}{}{\normalfont}{}{\bfseries}{}{ }{}
\theoremstyle{prim}

\newtheoremstyle{stil}{}{}{\slshape}{}{\bfseries}{}{ }{}
\theoremstyle{stil}
\newtheorem{thm}{Theorem}[section]
\newtheoremstyle{defi}{}{}{}{}{\bfseries}{}{ }{}
\theoremstyle{defi}
\newtheorem{defn}[thm]{Definition}
\theoremstyle{defi}
\newtheorem{rem}[thm]{Remark}
\theoremstyle{stil}
\newtheorem{pro}[thm]{Proposition}
\theoremstyle{stil}
\newtheorem{lem}[thm]{Lemma}
\theoremstyle{stil}
\newtheorem{kor}[thm]{Corollary}
\newenvironment{dok}{\noindent \textit{Proof.}}{\null\hfill$\qed$\hskip 2mm\vskip 2mm}

\newcommand{\rez}{\mathop{\mathrm{Res}}}
\newcommand{\om}{\mathop{\mathrm{Hom}}}
\newcommand{\ndo}{\mathop{\mathrm{End}}}
\newcommand{\sym}{\mathop{\mathrm{Sym}}}
\newcommand{\gauss}[2]{\genfrac{[}{]}{0pt}{}{#1}{#2}}
\newcommand{\diag}{\mathop{\mathrm{diag}}}

\numberwithin{equation}{section}

\begin{document}

\title[Principal subspaces for quantum affine algebra $U_q (A_n^{(1)})$]{Principal subspaces for quantum affine algebra $U_q(A_n^{(1)})$}

\author{Slaven Ko\v{z}i\'{c}}

\address{Department of Mathematics, University of Zagreb, 10000 Zagreb, Croatia}

\email{kslaven@math.hr}



\begin{abstract} 
We consider principal subspace $W(\Lambda)$ of integrable highest weight module $L(\Lambda)$ for quantum affine
algebra $U_q(\widehat{\mathfrak{sl}}_{n+1})$. We introduce quantum analogues of the quasi-particles 
associated with the principal subspaces for $\widehat{\mathfrak{sl}}_{n+1}$
 and discover certain relations
among them. By using these relations we find, for certain highest weight $\Lambda$, combinatorial bases of principal subspace
$W(\Lambda)$ 
in terms of monomials of quantum quasi-particles.
\end{abstract}

\maketitle


\section*{Introduction}

In \cite{D} Drinfeld discovered a remarkable realization of quantum affine algebras  
in the form of current operators.
In this paper we will study such a realization of  a quantum affine algebra $U_{q}(\widehat{\mathfrak{sl}}_{n+1})$.
Although being a Hopf algebra, Drinfeld proposed a new comultiplication formula (\cite{DF3}, \cite{DI}) for its new realization. 
This comultiplication formula proves to be very simple, when expressed in terms of the current operators, 
while the comultiplication formula induced from the usual Hopf algebra structure does not have a closed form in terms of those operators. 
In \cite{FJ} Frenkel and Jing found a realization of level $1$ integrable highest weight modules for quantum affine algebras of type $(ADE)^{(1)}$.
Similar realizations were obtained for integrable highest weight modules of some other algebra types, 
see for example \cite{I} and \cite{J} (level $2$ modules for $U_{q}(A_{1}^{(1)})$), 
\cite{B} and \cite{JM} (level $1$ modules for $U_{q}(B_{n}^{(1)})$), 
\cite{JKM} and \cite{JKM2} (level $-1/2$ and $1$ modules for $U_{q}(C_{n}^{(1)})$).
In \cite{DM} Ding and Miwa discovered a condition of quantum integrability, another important result formulated in terms of current operators.
This is actually an analogue of integrability relations for Lie algebra $\widehat{\mathfrak{sl}}_{2}$, found by Lepowsky and Primc (\cite{LP}) that 
 states that on every level $c$ integrable module of $U_{q}(\widehat{\mathfrak{sl}}_{n+1})$ we have 
$x_{\alpha_{i}}^{+}(z_{1})x_{\alpha_{i}}^{+}(z_{2})\cdots x_{\alpha_{i}}^{+}(z_{c+1})=0$ if 
$z_{1}/z_{2}=\ldots=z_{c}/z_{c+1}=q^{-2}$. Ding and Feigin formulated a similar result for their commutative quantum current operators
(\cite{DF}). Most of the above-mentioned results for algebra $U_{q}(A_{n}^{(1)})$ will play an important role in proving the main result 
of this paper.

One fundamental problem concerning quantum affine algebras is to associate with them a certain quantum vertex algebra theory
in the similar way vertex algebras are associated with affine Lie algebras. Several quantum vertex algebra theories
were developed: \cite{AB}, \cite{B2}, \cite{EK}, \cite{FR}, \cite{Li} and \cite{Li2}.
At the end of this paper we will use quantum vertex algebra 
theory discovered by Li (\cite{Li}, \cite{Li2}) in order to generalize certain objects studied in the preceding sections.

We will now briefly describe the main result of this paper. Let $\mathfrak{g}$ be a simple Lie algebra with a triangular decomposition
$$\mathfrak{g}=\mathfrak{n}_{-}\otimes \mathfrak{h}\otimes\mathfrak{n}_{+},$$ where $\mathfrak{h}$ is its Cartan subalgebra and
$\mathfrak{n}_{+}$ ($\mathfrak{n}_{-}$) is a direct sum of its one dimensional subalgebras corresponding to its positive (negative) roots. 
Denote by $\hat{\mathfrak{g}}$ an (untwisted) affine Lie algebra 
$\hat{\mathfrak{g}}=\mathfrak{g}\otimes\mathbb{C}[t,t^{-1}]\oplus \mathbb{C}c\oplus\mathbb{C}d$ and consider its 
decomposition $$\hat{\mathfrak{g}}=\hat{\mathfrak{n}}_{-}\otimes \hat{\mathfrak{h}}\otimes\hat{\mathfrak{n}}_{+},$$ where
 $\hat{\mathfrak{h}}=\mathfrak{h}\otimes\mathbb{C}[t,t^{-1}]\oplus \mathbb{C}c\oplus\mathbb{C}d$
and 
$\hat{\mathfrak{n}}_{\pm}=\mathfrak{n}_{\pm}\otimes\mathbb{C}[t,t^{-1}]$.
For a dominant integral highest weight $\Lambda$ denote by $L(\Lambda)$ an integrable highest weight module of $\hat{\mathfrak{g}}$
whose highest weight equals $\Lambda$. Define a principal subspace $W(\Lambda)$, 
 $$W(\Lambda):=U(\hat{\mathfrak{n}}_{+})\cdot v_{\Lambda},$$
 where $v_{\Lambda}$ is a highest weight vector of $L(\Lambda)$. Principal subspaces for an affine Lie algebra of type $A_{1}^{(1)}$
 were introduced by Feigin and Stoyanovsky (\cite{FS}). The authors constructed bases for the principal subspaces of $L(\Lambda)$
 consisting of the vectors 
 of the form $bv_{\Lambda}$, where $b$ is a monomial of quasi-particles. Georgiev (\cite{G}) generalized their results on
 affine Lie algebras of type $A_{n}^{(1)}$, $n\geq 1$, by constructing bases of the same type for the principal subspaces of 
 $L(\Lambda)$ for dominant integral weights $\Lambda$ of type
 \begin{equation}\label{000_intro}
\Lambda=c_{0}\Lambda_{0}+c_{j}\Lambda_{j},
\end{equation} 
where $j=1,2,\ldots,n$, $c_{0},c_{j}\in\mathbb{Z}_{\geq 0}$. 
The algebra $U_{q}(\hat{\mathfrak{g}})$ admits a decomposition
\begin{equation*}
U_{q}(\hat{\mathfrak{g}})\cong U_{q}(\hat{\mathfrak{n}}_{-})\otimes U_{q}(\hat{\mathfrak{h}})_{0}\otimes U_{q}(\hat{\mathfrak{n}}_{+})\quad\textrm{(vector space isomorphism)},
\end{equation*}
where $U_{q}(\hat{\mathfrak{n}}_{\pm})$ and $U_{q}(\hat{\mathfrak{h}})_{0}$ are subalgebras generated by certain elements given by Drinfeld realization.
Abusing the notation we denote by $L(\Lambda)$ an integrable highest weight module of 
$U_{q}(\hat{\mathfrak{g}})$, 
whose highest weight equals $\Lambda$ (and whose character is equal to the character of $\hat{\mathfrak{g}}$-module $L(\Lambda)$ (\cite{L3})).
Define a principal subspace $W(\Lambda)$ of the integrable highest weight module $L(\Lambda)$ of $U_{q}(\hat{\mathfrak{g}})$,
 $$W(\Lambda):=U_{q}(\hat{\mathfrak{n}}_{+})\cdot v_{\Lambda},$$
 where $v_{\Lambda}$ is a highest weight vector of $L(\Lambda)$. 
 The main goal of this paper is  to construct a monomial basis for $W(\Lambda)$ for quantum affine algebra 
 $U_{q}(\widehat{\mathfrak{sl}}_{n+1})$, where $\Lambda$ is an integral dominant highest weight as in (\ref{000_intro}).
 Such bases will be given in the following theorem: 
 
 \vspace{9pt}
 \noindent\textbf{Theorem \ref{main}}\textsl{ For a given highest weight $\Lambda=c_{0}\Lambda_{0}+c_{j}\Lambda_{j}$ as in (\ref{000_intro}) the sets 
$$\left\{bv_{\Lambda}\hspace{2pt}\left|\hspace{2pt}\right. b\in\mathfrak{B}_{W(\Lambda)}\right\}\qquad\textrm{and} \qquad \left\{bv_{\Lambda}\hspace{2pt}\left|\hspace{2pt}\right. b\in\bar{\mathfrak{B}}_{W(\Lambda)}\right\}$$ 
form the bases for the principal subspace $W(\Lambda)$. }
 \vspace{9pt}

These sets, $\mathfrak{B}_{W(\Lambda)}$ and $\bar{\mathfrak{B}}_{W(\Lambda)}$ consist of monomials of type $1$ or type $2$ quasi-particles whose charges and degrees
satisfy   certain difference conditions (see (\ref{t2_basis}) i (\ref{t1_basis})).

This paper is organized as follows. In Section 1 we recall some fundamental results concerning quantum affine 
algebras and their representations.

In Section 2 we define a principal subspace as above  and  study
its spanning sets.
Next, we introduce a notion of a quasi-particle of type $1$ and $2$.
Both of the types can be considered as a quantum analogue of quasi-particles studied in \cite{G}. 
The quasi-particles of type $1$ will be defined with the use of 
Frenkel-Jing realization and 
they will
allow us to carry out the linear independence part of the proof of Theorem 
\ref{main}.
Nevertheless, they will lack an important property: type $1$ quasi-particles 
of the same color do not commute. This fact will be 
our main motivation for introducing the quasi-particles of type $2$. 
They will be defined by using Ding-Feigin commutative
current operators (\cite{DF}) and they will commute when being of the same color. 
This property will prove to be crucial in the next section.

In Section 3 we study  relations among quasi-particles and then use these relations  to construct 
a spanning set $\bar{\mathfrak{B}}_{W(\Lambda)}$ of the principal subspace $W(\Lambda)$ (Theorem \ref{spanning_set}). The proof will be carried out in a similar way as in \cite{G}.  
Altogether there will be three sets of relations.
First, we will formulate the condition of quantum integrability (\cite{DM}) in terms of the quasi-particles, thus getting a first set of relations.
Next, we will discover relations among quasi-particles of adjacent colors. Some of the intermediate results will concern the type $1$
quasi-particles, but the final result will be formulated  in terms of type $2$ quasi-particles only. Finally, we will find the relations
among type $2$ quasi-particles of the same color.

Section 4 is devoted to the proof of   linear independence of the set $\mathfrak{B}_{W(\Lambda)}$ (Theorem \ref{independent}). 
Our main tools will be a certain projection $\pi_{(r_{n}^{(1)},\ldots,r_{1}^{(c)})}$
of the principal subspace and an operator $\mathcal{Y}(e^{\lambda_{i}},z)$ that can be considered as a certain modification
of the level $1$ vertex operators  realized by Koyama in \cite{Koyama}. The proof will be carried out as in \cite{G}.

Our definition of quasi-particles of type $1$, on an integrable highest weight module, heavily relies upon 
 Frenkel-Jing realization. 
However, we will show in Section 5 that type $1$ quasi-particles can be generalized to an arbitrary restricted module.
This generalization will arise as a corollary of quantum vertex algebra theory discovered by Li (\cite{Li}, \cite{Li2}).


\section{Preliminaries}

\subsection{Quantum affine algebra \texorpdfstring{$U_{q}(\hat{\mathfrak{g}})$}{Uq(g)}} 
We recall some facts from the theory of affine Kac-Moody Lie algebras (see \cite{Kac} for more details). 
Let $\hat{A}=(a_{ij})_{i,j=0}^{n}$ be a generalized Cartan matrix of affine type and let $S=\diag(s_{0},s_1,\ldots,s_n)$ be a diagonal matrix
of relatively prime positive integers such that the matrix $S\hat{A}$ is symmetric.

Let $\hat{\mathfrak{h}}$ be a vector space over $\mathbb{C}(q^{1/2})$
with a base $\left\{\alpha^{\vee}_{0},\alpha^{\vee}_{1},\ldots,\alpha^{\vee}_{n},d\right\}$. Denote by $\alpha_{0}$, $\alpha_{1}$, \ldots, $\alpha_{n}$ linear functionals from $\hat{\mathfrak{h}}^{*}$  such that
$$\alpha_{i}(\alpha_{j}^{\vee})=a_{ji},\quad\alpha_{i}(d)=\delta_{i0}$$
for $i,j=0,1\ldots,n$. 
Define a set of simple roots $\hat{\Pi}$ and a set of simple coroots  $\hat{\Pi}^{\vee}$ as
$$\hat{\Pi}:=\left\{\alpha_{i}\hspace{2pt}|\hspace{2pt}i=0,1,\ldots,n\right\}\qquad\textrm{and}\qquad \hat{\Pi}^{\vee}:=\left\{\alpha_{i}^{\vee}\hspace{2pt}|\hspace{2pt}i=0,1,\ldots,n\right\}.$$
Then the ordered triple $(\hat{\mathfrak{h}},\hat{\Pi},\hat{\Pi}^{\vee})$ is a realization of the matrix $\hat{A}$. We will denote by 
$\hat{\mathfrak{g}}$
the affine Kac-Moody Lie algebra associated with the matrix $\hat{A}$.

Define a root lattice $\hat{Q}$ and a set  $\hat{Q}^{+}$ as 
$$\hat{Q}:=\bigoplus_{i=0}^{n}\mathbb{Z}\alpha_{i}\quad\textrm{and}\quad \hat{Q}^{+}:=\bigoplus_{i=0}^{n}\mathbb{Z}_{\geq 0}\alpha_{i}.$$ 
Denote by $\Lambda_{0}$, $\Lambda_{0}$, \ldots, $\Lambda_{n}$ elements of $\hat{\mathfrak{h}}^{*}$ such that
$$\Lambda_{i}(\alpha_{j}^{\vee})=\delta_{ij},\quad \Lambda_{i}(d)=0$$
for $i,j=0,1,\ldots,n$. The center of the Lie algebra $\hat{\mathfrak{g}}$ is one-dimensional and it is generated by the element
$$c=c_{0}\alpha_{0}^{\vee}+c_{1}\alpha_{1}^{\vee}+\ldots+c_{n}\alpha_{n}^{\vee}\in \hat{\mathfrak{h}}$$
and imaginary roots of $\hat{\mathfrak{g}}$ are integer multiples of
$$\delta=d_{0}\alpha_{0}+d_{1}\alpha_{1}+\ldots+d_{n}\alpha_{n}\in \hat{\mathfrak{h}}^{*},$$
where integers $c_i$ and $d_i$ are given in \cite{Kac}.
If $\hat{\mathfrak{g}}=\widehat{\mathfrak{sl}}_{n+1}$ we have $c_{i}=d_{i}=1$.
Define a weight lattice $\hat{P}$ as a free Abelian group generated by the elements $\Lambda_{0}$, $\Lambda_{1}$, \ldots, $\Lambda_{n}$ and $\frac{1}{d_{0}}\delta$.
Define an integral dominant weight as the element of $\hat{P}^{+}$,
$$\hat{P}^{+}:=\left\{\Lambda\in \hat{P}\hspace{2pt}|\hspace{2pt}\Lambda(\alpha_{i}^{\vee})\in\mathbb{Z}_{\geq 0}\textrm{ for }i=0,1,\ldots,n\right\}.$$
The invariant symmetric bilinear form on $\hat{\mathfrak{h}}^{*}$ is given by
$$(\alpha_{i},\alpha_{j})=s_{i}a_{ij},\quad (\delta,\alpha_{i})=(\delta,\delta)=0$$
for $i,j=0,1,\ldots,n$.

Denote by $\mathfrak{g}$ a simple Lie algebra associated with the Cartan matrix 
$A=(a_{ij})_{i,j=1}^{n}$.  Let $\mathfrak{h}\subset\hat{\mathfrak{h}}$ be a Cartan subalgebra of  
$\mathfrak{g}$, which is generated by the elements $\alpha_{1}^{\vee}$, $\alpha_{2}^{\vee}$, \ldots, $\alpha_{n}^{\vee}$.
Denote by 
$$Q:=\bigoplus_{i=1}^{n}\mathbb{Z}\alpha_{i}\subset\mathfrak{h}^{*}\quad\textrm{and}\quad P:=\bigoplus_{i=1}^{n}\mathbb{Z}\lambda_{i}\subset\mathfrak{h}^{*}$$ 
the classical root lattice and the classical weight lattice (associated with $\mathfrak{g}$), where elements $\lambda_{i}$ 
satisfy $\lambda_{i}(\alpha_{j}^{\vee})=\delta_{ij}$
for $i,j=1,2,\ldots,n$.

Fix an indeterminate $q$.
For any two integers $m$ and $k$, $k>0$,   define $q$-integers,
 $$[m]_q:=\frac{q^{m}-q^{-m}}{q-q^{-1}}, $$
 and $q$-factorials,
 $$[0]_q !:=1,\quad [k]_q !:=[k]_{q}[k-1]_{q}\cdots[1]_{q}.$$
 For all nonnegative integers $m$ and $k$, $m\geq k$, define $q$-binomial coefficients as
 $$\gauss{m}{k}_{q}:=\frac{[m]_{q}!}{[k]_{q}![m-k]_{q}!}.$$

We  present a Drinfeld realization of the quantum affine algebra $U_{q}(\hat{\mathfrak{g}})$.

\begin{defn}[\cite{D}]\label{drinfeld}
The quantum affine algebra $U_{q}(\hat{\mathfrak{g}})$ is the associative algebra over $\mathbb{C}(q^{1/2})$ with unit $1$ generated by the 
elements
$x_{\alpha_{i}}^{\pm}(k)$, $a_{i}(l)$, $K_{i}^{\pm 1}$, $\gamma^{\pm 1/2}$ and $q^{\pm d}$, $i=1,2,\ldots,n$, $k,l\in\mathbb{Z}$, $l\neq 0$, subject to the following relations:
\begin{align*}
& [\gamma^{\pm1/2},u]=0\textrm{ for all }u\in U_{q}(\hat{\mathfrak{g}})_0,\\
& K_i K_j=K_j K_i,\quad K_i K_{i}^{-1}=K_{i}^{-1}K_{i}=1,\\
& [a_{i}(k),a_{j}(l)]=\delta_{k+l\hspace{2pt}0}\frac{[a_{ij}k]_{q_{i}}}{k}\frac{\gamma^{k}-\gamma^{-k}}{q_{j}-q_{j}^{-1}},\\
& [a_{i}(k),K_{j}^{\pm 1}]=[q^{\pm d},K_{j}^{\pm 1}]=0,\\
& q^{d}x_{\alpha_i}^{\pm}(k)q^{-d}=q^{k}x_{\alpha_i}^{\pm}(k),\quad q^{d}a_{i}(l)q^{-d}=q^{k}a_{i}(l),\\
& K_{i}x_{\alpha_j}^{\pm}(k)K_{i}^{-1}=q^{\pm (\alpha_{i},\alpha_{j})}x_{\alpha_j}^{\pm }(k) ,\\
& [a_{i}(k),x_{\alpha_{j}}^{\pm}(l)]=\pm\frac{[a_{ij}k]_{q_{i}}}{k}\gamma^{\mp |k|/2}x_{\alpha_{j}}^{\pm}(k+l),\\
& x_{\alpha_{i}}^{\pm}(k+1)x_{\alpha_{j}}^{\pm}(l)-q^{\pm(\alpha_i,\alpha_j)}x_{\alpha_{j}}^{\pm}(l)x_{\alpha_{i}}^{\pm}(k+1)\\
&\hspace{20pt}=q^{\pm(\alpha_i,\alpha_j)}x_{\alpha_{i}}^{\pm}(k)x_{\alpha_{j}}^{\pm}(l+1)-x_{\alpha_{j}}^{\pm}(l+1)x_{\alpha_{i}}^{\pm}(k),\\
& [x_{\alpha_i}^{+}(k),x_{\alpha_{j}}^{-}(l)]=\frac{\delta_{ij}}{q_{i}-q^{-1}_{i}}\left(\gamma^{\frac{k-l}{2}}\psi_{i}(k+l)-\gamma^{\frac{l-k}{2}}\phi_{i}(k+l)\right),\\
& \sym_{l_1,l_2,\ldots,l_m}\sum_{s=0}^{m}(-1)^{s}\gauss{m}{s}_{q_i}x_{\alpha_{i}}^{\pm}(l_{1})\cdots x_{\alpha_{i}}^{\pm}(l_{s})x_{\alpha_{j}}^{\pm}(k)x_{\alpha_{i}}^{\pm}(l_{s+1})\cdots x_{\alpha_{i}}^{\pm}(l_{m})=0, \quad\textrm{for }i\neq j,
\end{align*}
where $m:=1-a_{ij}$.
The elements $\phi_{i}(-r)$ and $\psi_{i}(r)$, $r\in\mathbb{Z}_{\geq 0}$, are given by \\
\begin{align*}
& \phi_{i}(z):=\sum_{r=0}^{\infty}\phi_{i}(-r)z^{r}:=K_{i}^{-1}\exp\left(-(q_{i}-q_{i}^{-1})\sum_{r=1}^{\infty}a_{i}(-r)z^{r}\right),\\
& \psi_{i}(z):=\sum_{r=0}^{\infty}\psi_{i}(r)z^{-r}:=K_{i}\exp\left((q_{i}-q_{i}^{-1})\sum_{r=1}^{\infty}a_{i}(r)z^{-r}\right).
\end{align*}
\end{defn}

If $q_{i}=q$ we will usually omit the index $q_i$ and write $[m]$ instead of $[m]_{q_{i}}$. Denote by $x_{\alpha_{i}}^{\pm}(z)$ the series
\begin{equation}\label{101_exp:series}
x_{\alpha_{i}}^{\pm}(z)=\sum_{r\in\mathbb{Z}}x_{\alpha_{i}}^{\pm}(r)z^{-r-1}\in U_{q}(\hat{\mathfrak{g}})[[z^{\pm 1}]].
\end{equation}
We shall continue to use the notation $x_{\alpha_{i}}^{\pm}(z)$ for the action of the expression (\ref{101_exp:series}) on an arbitrary
$U_{q}(\hat{\mathfrak{g}})$-module $V$:
$$x_{\alpha_{i}}^{\pm}(z)=\sum_{r\in\mathbb{Z}}x_{\alpha_{i}}^{\pm}(r)z^{-r-1}\in (\ndo V)[[z^{\pm 1}]].$$

Let $i=1,2,\ldots, n$. Denote by $U_{q}(\hat{\mathfrak{n}}^{\pm}_{i})$ a subalgebra of $U_{q}(\hat{\mathfrak{g}})$ generated by the elements $x^{\pm}_{\alpha_{i}}(m)$, $m\in\mathbb{Z}$.
Denote by $U_{q}(\hat{\mathfrak{n}}^{\pm})$ a subalgebra of $U_{q}(\hat{\mathfrak{g}})$ generated by the elements $x^{\pm}_{\alpha_{i}}(m)$, $m\in\mathbb{Z}$, $i=1,2,\ldots, n$.
Finally, denote by $U_{q}(\hat{\mathfrak{h}})_{0}$ a subalgebra of $U_{q}(\hat{\mathfrak{g}})$ generated by the elements $a_{i}(l)$, $K_{i}^{\pm 1}$, $\gamma^{\pm 1/2}$ and $q^{\pm d}$ for $i=1,2,\ldots,n$, $l\in\mathbb{Z}$, $l\neq 0$.
Multiplication establishes an isomorphism of $\mathbb{C}(q^{1/2})$-vector spaces:
\begin{equation}\label{102_triangle}
U_{q}(\hat{\mathfrak{g}})\cong U_{q}(\hat{\mathfrak{n}}^{-})\otimes U_{q}(\hat{\mathfrak{h}})_{0}\otimes U_{q}(\hat{\mathfrak{n}}^{+}).
\end{equation}

Drinfeld gave the Hopf algebra structure for the realization \ref{drinfeld} of the algebra  $U_{q}(\widehat{\mathfrak{sl}}_{n+1})$.

\begin{thm}
The algebra $U_{q}(\widehat{\mathfrak{sl}}_{n+1})$ has a Hopf algebra structure which is given by the following formulas for the 

coproduct $\Delta$:
\begin{align*}
&\Delta(q^{c/2})=q^{c/2}\otimes q^{c/2},\\
&\Delta(x_{\alpha_{i}}^{+}(z))=x_{\alpha_{i}}^{+}(z)\otimes 1 + \phi_{i}(zq^{c_{1}/2})\otimes x_{\alpha_{i}}^{+}(zq^{c_{1}}),\\
&\Delta(x_{\alpha_{i}}^{-}(z))=1\otimes x_{\alpha_{i}}^{-}(z)+x_{\alpha_{i}}^{-}(zq^{c_{2}})\otimes\psi_{i}(zq^{c_{2}/2}),\\
&\Delta(\phi_{i}(z))=\phi_{i}(zq^{-c_{2}/2})\otimes \phi_{i}(zq^{c_{1}/2}),\\
&\Delta(\psi_{i}(z))=\psi_{i}(zq^{c_{2}/2})\otimes \psi_{i}(zq^{-c_{1}/2}),
\end{align*}
where $q^{c_{1}}$ means the action of the center $q^{c}$ on the first tensor component and $q^{c_{2}}$ means the action of the center $q^{c}$ on the second component;

counit $\varepsilon$:
\begin{align*}
&\varepsilon(q^c)=1, \quad \varepsilon(x^{\pm}_{\alpha_{i}}(z))=0,\quad\varepsilon(\phi_{i}(z))=\varepsilon(\psi_{i}(z))=1;
\end{align*}

antipode $S$:
\begin{align*}
& S(q^{c})=q^{-c},\\
& S(x_{\alpha_{i}}^{+}(z))=-\phi_{i}(zq^{-c/2})^{-1}x_{\alpha_{i}}^{+}(zq^{-c}),\\
& S(x_{\alpha_{i}}^{-}(z))=-x_{\alpha_{i}}^{-}(zq^{-c})\psi_{i}(zq^{-c/2})^{-1},\\
& S(\phi_{i}(z))=\phi_{i}(z)^{-1},\quad S(\psi_{i}(z))=\psi_{i}(z)^{-1}.
\end{align*}
\end{thm}

The comultiplication structure requires a certain completion of the tensor space. However, the comultiplication  is well defined
if at least one of the tensor factors is a highest weight module. The proof for the the above theorem for the case of 
$U_{q}(\widehat{\mathfrak{sl}}_{2})$ is given in \cite{DI}.
For $l\in\mathbb{Z}_{>0}$ set
\begin{equation*}
\Delta^{(0)}:=1\quad\textrm{ and }\quad\Delta^{(l)}:=(\underbrace{1\otimes\cdots\otimes 1}_{\mbox{$l-1$}} \otimes \Delta)\Delta^{(l-1)}.
\end{equation*}
The coproduct formula applied on the tensor product of $c$  highest weight modules gives
\begin{equation}
\label{103_coproduct}
\Delta^{(c-1)}(x_{\alpha_{i}}^{+}(z))=\sum_{l=1}^{c}x_{\alpha_{i}}^{+(l;c)}(z),
\end{equation}
where
\begin{equation*}
x_{\alpha_{i}}^{+(l;c)}(z)=\underbrace{\phi_{i}(zq^{\frac{1}{2}})\otimes\phi_{i}(zq^{\frac{3}{2}})\otimes\cdots\otimes\phi_{i}(zq^{l-\frac{3}{2}})}_{\mbox{$l-1$}}\otimes x_{\alpha_{i}}^{+}(zq^{l-1})\otimes\underbrace{1\otimes\cdots\otimes 1}_{\mbox{$c-l$}}.
\end{equation*}


\subsection{Frenkel-Jing construction}
Let $V$ be a $U_{q}(\hat{\mathfrak{g}})$-module. The module $V$ is called a level $c$ module if the central element $\gamma$ acts as the scalar $q^c$ 
on $V$. We will study the integrable highest weight modules $L(\Lambda)$ associated with the integral dominant highest weights $\Lambda\in\hat{P}^{+}$.
The level of such a module $L(\Lambda)$ is a positive integer. Some basic facts about this modules can be found in \cite{HK}.
From now on we will assume that the Lie algebra $\hat{\mathfrak{g}}$ is of type $A_{n}^{(1)}$. We present the Frenkel-Jing construction
of level $1$ integrable highest weight modules (\cite{FJ}). 

Let $V$ be a $U_{q}(\hat{\mathfrak{g}})$-module of level $c$. The Heisenberg algebra  $U_{q}(\hat{\mathfrak{h}})$ of level $c$ is generated by the elements
 $a_{i}(k)$, $i=1,2,\ldots,n$, $k\in\mathbb{Z}\setminus\left\{0\right\}$ and the central element $\gamma^{\pm 1}=q^{\pm c}$ subject to the relations 
\begin{equation}\label{104_heisenberg}
[a_{i}(r),a_{j}(s)]=\delta_{r+s\hspace{2pt}0}\frac{[a_{ij}r][cr]}{r}
\end{equation}
for $i,j=1,2,\ldots,n$, $r,s\in\mathbb{Z}\setminus\left\{0\right\}$. $V$ is also a $U_{q}(\hat{\mathfrak{h}})$-module.

 We define the following operators on $V$:
 \begin{align*}
 &E_{-}^{\pm}(a_{i},z):=\exp\left(\mp\sum_{r\geq 1}\frac{q^{\mp cr/2}}{[cr]}a_{i}(-r)z^{r}\right),\\
  &E_{+}^{\pm}(a_{i},z):=\exp\left(\pm\sum_{r\geq 1}\frac{q^{\mp cr/2}}{[cr]}a_{i}(r)z^{-r}\right).
\end{align*}

The Heisenberg algebra $U_{q}(\hat{\mathfrak{h}})$ has a natural realization on the space $\sym(\hat{\mathfrak{h}}^{-})$
of the symmetric algebra generated by the elements $a_{i}(-r)$, 
  $r\in\mathbb{Z}_{>0}$, $i=1,2,\ldots,n$ via the following rule:
  \begin{align*}
\gamma^{\pm 1}\hspace{5pt}&\ldots\hspace{5pt}\textrm{multiplication by }q^{\pm c},\\
a_{i}(r)\hspace{5pt}&\ldots\hspace{5pt}\textrm{differentiation operator subject to (\ref{104_heisenberg})},\\
a_{i}(-r)\hspace{5pt}&\ldots\hspace{5pt}\textrm{multiplication by the element }a_{i}(-r).
\end{align*}
We denote the resulted irreducible $U_{q}(\hat{\mathfrak{h}})$-module
 as $K(c)$.

The associative algebra $\mathbb{C}\left\{P\right\}$ (of the classical weight lattice $P$)
is generated by the elements $e^{\alpha_2},e^{\alpha_{3}},\ldots e^{\alpha_{n}}$ and $e^{\lambda_{n}}$ subject to the relations
\begin{align*}
&e^{\alpha_{i}}e^{\alpha_{j}}=(-1)^{(\alpha_i,\alpha_j)}e^{\alpha_{j}}e^{\alpha_{i}},\\
&e^{\alpha_{i}}e^{\lambda_{n}}=(-1)^{\delta_{in}}e^{\lambda_{n}}e^{\alpha_{i}},
\end{align*}
$i,j=2,3,\ldots,n$. For 
$\alpha=m_{2}\alpha_2+ m_{3}\alpha_3+\ldots+m_{n}\alpha_{n}+m_{n+1}\lambda_{n}\in P$
we denote $e^{m_{2}\alpha_2}e^{m_{3}\alpha_3}\cdots e^{m_{n}\alpha_{n}}e^{m_{n+1}\lambda_{n}}\in \mathbb{C}\left\{P\right\}$ by
 $e^{\alpha}$. For example, we have
\begin{align*}
& e^{\alpha_{1}}:=e^{-2\alpha_{2}}e^{-3\alpha_{3}}\cdots e^{-n\alpha_{n}}e^{(n+1)\lambda_{n}},\\
& e^{\lambda_{i}}:=e^{-\alpha_{i+1}}e^{-2\alpha_{i+2}}\cdots e^{-(n-i)\alpha_{n}}e^{(n+1)\lambda_{n}}
\end{align*}
for $i=0,1,\ldots,n-1$. 

One can easily verify the following relations:
\begin{pro} For $i,j=1,2,\ldots,n$ we have
\begin{enumerate}
\item $e^{\alpha_{1}}e^{\lambda_{n}}=(-1)^{n}e^{\lambda_{n}}e^{\alpha_{1}}$,
\item $e^{\alpha_{i}}e^{\alpha_{j}}=(-1)^{(\alpha_i,\alpha_j)}e^{\alpha_{j}}e^{\alpha_{i}}.$\\
\end{enumerate}
\end{pro}

Denote by $\mathbb{C}\left\{Q\right\}$ the subalgebra of $\mathbb{C}\left\{P\right\}$ generated by the elements $e^{\alpha_i}$, $i=1,2,\ldots,n$.
Define the vector space
$$\mathbb{C}\left\{Q\right\}e^{\lambda_i}:=\left\{ae^{\lambda_{i}}\hspace{2pt}|\hspace{2pt}a\in\mathbb{C}\left\{Q\right\}\right\}.$$
For $\alpha\in Q$ define an action  $z^{\partial_\alpha}$ on $\mathbb{C}\left\{Q\right\}e^{\lambda_i}$ by
$$z^{\partial_\alpha}e^{\beta}e^{\lambda_i}=z^{(\alpha,\beta+\lambda_i)}e^{\beta}e^{\lambda_i}.$$

\begin{thm}[\cite{FJ}]
 By the action 
\begin{align*}
x_{\alpha_j}^{\pm}(z)&:=E_{-}^{\pm}(-a_{j},z)E_{+}^{\pm}(-a_{j},z)\otimes e^{\pm\alpha_{j}}z^{\pm\partial_{\alpha_{j}}},
\end{align*}
$j=1,2,\ldots,n$, 
the space $$L_{i}:=K(1)\otimes\mathbb{C}\left\{Q\right\}e^{\lambda_i}$$ becomes
the integrable highest weight module of $U_{q}(\widehat{\mathfrak{sl}}_{n+1})$  with the highest weight $\Lambda_i$.
The vector $1\otimes e^{\lambda_{i}}$ is the highest weight vector of $L_{i}$.
\end{thm}


\subsection{Ding-Feigin operators}
We briefly sketch the construction of the commutative quantum current operators $\bar{x}_{\alpha_i}^{+}(z)$ from \cite{DF}.
First, we define the following operators on any level $c$ integrable highest weight module of $U_{q}(\widehat{\mathfrak{sl}}_{n+1})$:
\begin{equation*}
k^{+}_{i}(z):=\exp\left((q-q^{-1})\sum_{r\geq 1}\frac{-q^{2r+\frac{cr}{2}}}{1+q^{2r}}a_{i}(r)z^{-r}\right),
\end{equation*}
$i=1,2,\ldots,n$.
Next, we define operators
\begin{equation}\label{105_ding_feigin}
\bar{x}_{\alpha_i}^{+}(z):=x_{\alpha_i}^{+}(z)k^{+}_{i}(z).
\end{equation}

\begin{thm}[\cite{DF}]\label{106_ding_feigin}
For $i,j=1,2,\ldots,n$ such that $a_{ij}=0,2$
we have
\begin{align*}
&\bar{x}_{\alpha_i}^{+}(z_1)\bar{x}_{\alpha_j}^{+}(z_2)=\bar{x}_{\alpha_j}^{+}(z_2)\bar{x}_{\alpha_i}^{+}(z_1).
\end{align*}
\end{thm}

At the end  we list a few formulas which will prove to be useful in next sections. 
Most of them can be found in \cite{DM} or \cite{DF} and also, they can be verified by a direct calculation.

\begin{lem} Let $i,j=1,2,\ldots,n$ such that $a_{ij}=-1$.
\begin{enumerate}
\item On every level $c$ integrable highest weight module we have
\begin{align}
&\hspace{-75pt}\displaystyle k_{i}^{+}(z_{1})x_{\alpha_i}^{+}(z_{2})=\frac{1-q^{2}\frac{z_{2}}{z_{1}}}{1-\frac{z_{2}}{z_{1}}}x_{\alpha_i}^{+}(z_{2})k_{i}^{+}(z_{1}),\label{f:1}\\
 &\hspace{-75pt} \displaystyle k_{i}^{+}(z_{1})x_{\alpha_j}^{+}(z_{2})=\exp\left(\sum_{r\geq 1}\frac{q^{r}-q^{-r}}{r(q^{2r}+1)}\left(\frac{z_{2}q^{2}}{z_{1}}\right)^{r}\right)x_{\alpha_j}^{+}(z_{2})k_{i}^{+}(z_{1}).\label{f:2}
\end{align}
\item On every level $1$ integrable highest weight module we have
\begin{align}
&\hspace{-75pt}\displaystyle x_{\alpha_i}^{+}(z_{1})x_{\alpha_i}^{+}(z_{2})=z_{1}^{2}\left(1-\frac{z_{2}}{z_{1}}\right)\left(1-q^{-2}\frac{z_{2}}{z_{1}}\right): x_{\alpha_i}^{+}(z_{1})x_{\alpha_i}^{+}(z_{2}):,\label{f:3}\\
  &\hspace{-75pt}\displaystyle x_{\alpha_i}^{+}(z_{1})\phi_{i}(z_{2}q^{1/2})=q^{2}\frac{1-q^{-2}\frac{z_{2}}{z_{1}}}{1-q^{2}\frac{z_{2}}{z_{1}}}\phi_{i}(z_{2}q^{1/2})x_{\alpha_i}^{+}(z_{1}),\label{f:4}\\
  &\hspace{-75pt}\displaystyle x_{\alpha_i}^{+}(z_{1})x_{\alpha_j}^{+}(z_{2})=\frac{1}{z_{1}}\frac{1}{1-q^{-1}\frac{z_{2}}{z_{1}}}: x_{\alpha_i}^{+}(z_{1})x_{\alpha_j}^{+}(z_{2}):,\label{f:5}\\
  &\hspace{-75pt}\displaystyle x_{\alpha_i}^{+}(z_{1})\phi_{j}(z_{2}q^{1/2})=q^{-1}\frac{1-q\frac{z_{2}}{z_{1}}}{1-q^{-1}\frac{z_{2}}{z_{1}}} \phi_{j}(z_{2}q^{1/2})x_{\alpha_i}^{+}(z_{1}).\label{f:6}
\end{align}
\end{enumerate}
\end{lem}


\section{Principal subspaces and quasi-particles}


\subsection{Principal subspaces}
We introduce the notion of a principal subspace.
Recall decomposition (\ref{102_triangle}).

\begin{defn}
Let $v_{\Lambda}$ be the highest weight vector of  $L(\Lambda)$, $\Lambda\in\hat{P}^{+}$. 
We define a principal subspace $W(\Lambda)$ of the integrable highest weight module $L(\Lambda)$ of $U_{q}(\widehat{\mathfrak{sl}}_{n+1})$,
$$W(\Lambda):=U_{q}({\mathfrak{\hat{n}^{+}}})v_{\Lambda}.$$
\end{defn}

\begin{lem}\label{colors}
For every $\Lambda\in\hat{P}^{+}$ we have
$$W(\Lambda)=U_{q}(\hat{\mathfrak{n}}^{+}_{n})U_{q}(\hat{\mathfrak{n}}^{+}_{n-1})\ldots U_{q}(\hat{\mathfrak{n}}^{+}_{1})v_{\Lambda}.$$
\end{lem}

\begin{dok}
It is sufficient to find a way to change the order in the product $x_{\alpha_{i}}^{+}(l_1)x_{\alpha_{j}}^{+}(l_2)$, $i,j=1,2,\ldots,n$, $l_1,l_2\in\mathbb{Z}$, when 
acting on a given vector $v\in W(\Lambda)$. Fix an integer $N$ such that $x_{\alpha_{j}}^{+}(l)v=0$ for $l>N$. Recall the formula
$$x_{\alpha_{i}}^{\pm}(k+1)x_{\alpha_{j}}^{\pm}(l)-q^{\pm(\alpha_i,\alpha_j)}x_{\alpha_{j}}^{\pm}(l)x_{\alpha_{i}}^{\pm}(k+1)
=q^{\pm(\alpha_i,\alpha_j)}x_{\alpha_{i}}^{\pm}(k)x_{\alpha_{j}}^{\pm}(l+1)-x_{\alpha_{j}}^{\pm}(l+1)x_{\alpha_{i}}^{\pm}(k)$$
from Definition \ref{drinfeld}. By applying the formula on a vector $v$  we get
\begin{align*}
&x^{+}_{\alpha_{i}}(l_1)x^{+}_{\alpha_{j}}(l_2)v\\
&\hspace{10pt}=q^{(\alpha_i,\alpha_j)}x^{+}_{\alpha_{j}}(l_2)x^{+}_{\alpha_{i}}(l_1)v-x^{+}_{\alpha_{j}}(l_2 +1)x^{+}_{\alpha_{i}}(l_1 -1)v+q^{(\alpha_i,\alpha_j)}x^{+}_{\alpha_{i}}(l_1 -1)x^{+}_{\alpha_{j}}(l_2 +1)v.
\end{align*}
The right hand side of the above equality contains two terms of the changed order, 
$x^{+}_{\alpha_{j}}(l_2)x^{+}_{\alpha_{i}}(l_1)v$ and
$x^{+}_{\alpha_{j}}(l_2 +1)x^{+}_{\alpha_{i}}(l_1 -1)v$,
 and one (last) term of the previous order,
$x^{+}_{\alpha_{i}}(l_1 -1)x^{+}_{\alpha_{j}}(l_2 +1)v$. 
Notice that the degree of  $x_{\alpha_j}^{+}$ in the last term equals $l_2 +1$ while the degree of $x_{\alpha_j}^{+}$ 
on the left hand side of the equality equals $l_2$. Therefore, by applying the same formula  $r$ more times, 
where $r\in\mathbb{Z}_{>0}$ is such that $l_{2}+r\geq N$, we will get the vector  $x^{+}_{\alpha_{i}}(l_1)x^{+}_{\alpha_{j}}(l_2)v$
written as a linear combination of some vectors $x^{+}_{\alpha_{j}}(l)x^{+}_{\alpha_{i}}(k)v$, thus 
finishing the proof.
\end{dok}

For $i=1,2,\ldots,n$ and a vector $v\in W(\Lambda)$  denote by $\bar{U}_{q}(\hat{\mathfrak{n}}^{+}_{i})v$
a subspace of $W(\Lambda)$ spanned by the vectors 
$\bar{x}^{+}_{\alpha_{i}}(l_1)\cdots \bar{x}^{+}_{\alpha_{i}}(l_k)v,$
where $k\in\mathbb{Z}_{\geq 0}$, $l_{1},\ldots,l_{k}\in\mathbb{Z}.$

\begin{lem}\label{colors_2}
For every $\Lambda\in\hat{P}^{+}$ we have
$$W(\Lambda)=\bar{U}_{q}(\hat{\mathfrak{n}}^{+}_{n})\bar{U}_{q}(\hat{\mathfrak{n}}^{+}_{n-1})\ldots \bar{U}_{q}(\hat{\mathfrak{n}}^{+}_{1})v_{\Lambda}.$$
\end{lem}

\begin{dok}
Every vector 
$$w:=\bar{x}_{\alpha_{n}}^{+}(l_{1,n})\ldots\bar{x}_{\alpha_{n}}^{+}(l_{k_{n},n})\bar{x}_{\alpha_{n-1}}^{+}(l_{1,n-1})\ldots \bar{x}_{\alpha_{1}}^{+}(l_{1,1})\ldots\bar{x}_{\alpha_{1}}^{+}(l_{k_{1},1})v_{\Lambda}$$ 
of the space $\bar{U}_{q}(\hat{\mathfrak{n}}^{+}_{n})\bar{U}_{q}(\hat{\mathfrak{n}}^{+}_{n-1})\ldots \bar{U}_{q}(\hat{\mathfrak{n}}^{+}_{1})v_{\Lambda}$
is a product of $k_n + \ldots +k_2 + k_1$ factors of the form
$\sum_{m\geq 0}x_{\alpha_{i}}^{+}(l-m)k_{i}^{+}(m)$, for some $i=1,2,\ldots , n$  and  $l\in\mathbb{Z}$,
acting on the highest weight vector $v_{\Lambda}$. Furthermore, every $k_{i}^{+}(m)$ equals certain linear combination of products of the elements
$a_{i}(s)$, $s\geq 1$. Recall $a_{i}(s)v_{\Lambda}=0$ for $s\geq 1$. Therefore, by using the formula 
\begin{equation*}\label{def2}
[a_{i}(k),x_{\alpha_{j}}^{+}(l)]=\frac{[a_{ij}k]}{k}\gamma^{-\frac{|k|}{2}}x_{\alpha_{j}}^{+}(k+l)
\end{equation*} 
(see Definition \ref{drinfeld}) we can write the vector $w$ as a linear combination of the vectors of the form
\begin{equation*}
x_{\alpha_{n}}^{+}(l^{'}_{1,n})\ldots x_{\alpha_{n}}^{+}(l^{'}_{k_{n},n})x_{\alpha_{n-1}}^{+}(l^{'}_{1,n-1})\ldots x_{\alpha_{1}}^{+}(l^{'}_{1,1})\ldots x_{\alpha_{1}}^{+}(l^{'}_{k_{1},1})v_{\Lambda}.
\end{equation*}
This implies
$\bar{U}_{q}(\hat{\mathfrak{n}}^{+}_{n})\bar{U}_{q}(\hat{\mathfrak{n}}^{+}_{n-1})\ldots \bar{U}_{q}(\hat{\mathfrak{n}}^{+}_{1})v_{\Lambda}\subseteq W(\Lambda)\textrm{.}$

The opposite inclusion can be proven in the same way by using the formula
$$[a_{i}(k),\bar{x}_{\alpha_{j}}^{+}(l)]=\frac{[a_{ij}k]}{k}\gamma^{-\frac{k}{2}}\bar{x}_{\alpha_{j}}^{+}(k+l)$$
for $i,j=1,2,\ldots  n$, $k,l\in\mathbb{Z}$, $k\geq 1$, and Lemma \ref{colors}. The above formula can be easily verified by a direct calculation. 
\end{dok}


\subsection{Quasi-particles of type 1}
For a given vector space $V$ and for any two Laurent series
$a(z)=\sum_{r\in\mathbb{Z}}a_{r}z^{-r-1}$ and $b(z)=\sum_{s\in\mathbb{Z}}b_{s}z^{-s-1}$ in $(\ndo V)[[z^{\pm 1}]]$ the product 
$a(z_1)b(z_2)$  is obviously well-defined
element of $(\ndo V)[[z_{1}^{\pm 1},z_{2}^{\pm 1}]]$. However,  the product $a(z)b(z)$  may not be a well-defined element of $(\ndo V)[[z^{\pm 1}]]$.
With an additional assumption of $a(z_1)b(z_2)$ being an element of $\om (V,V((z_1,z_2)))$ we can conclude that $a(z)b(z)\in\om(V,V((z)))\subset (\ndo V)[[z^{\pm 1}]]$.
More generally, if $a(z_1,\ldots,z_m)\in\om (V,V((z_1,\ldots, z_m)))$, we can carry out the substitution $z_{p}\to z$ for $p=1,2,\ldots,n$ and thus get
the element $a(z,\ldots,z)\in\om (V,V((z)))$. We will denote the application of such a substitution  by a limit symbol, for example
$$\lim_{z_{p}\to z}a(z_1,\ldots,z_m)=a(z,\ldots,z)\in\om (V,V((z))).$$ 

Fix $i=1,2,\ldots,n$.
On any integrable  level $1$ module $L(\Lambda_{j})$ of $U_{q}(\widehat{\mathfrak{sl}}_{n+1})$ we have
\begin{equation}\label{order}
\phi_{i}(z_1)\ldots \phi_{i}(z_r):x_{\alpha_{i}}^{+}(z_{r+1})\ldots x_{\alpha_{i}}^{+}(z_{m}):\hspace{2pt}\in \om(L(\Lambda_{j}),L(\Lambda_{j})((z_1,\ldots,z_m))).
\end{equation}
Therefore, there exists a limit $z_{p}\to z$, $p=1,2,\ldots,m$, of the above expression:
$$\lim_{z_{p}\to z}\left(\phi_{i}(z_1)\ldots \phi_{i}(z_r):x_{\alpha_{i}}^{+}(z_{r+1})\ldots x_{\alpha_{i}}^{+}(z_{m}):\right)\hspace{2pt}\in \om(L(\Lambda_{j}),L(\Lambda_{j})((z))).$$

Denote by $V$ a level $c$ integrable highest weight  module $L(\Lambda)$ of $U_{q}(\widehat{\mathfrak{sl}}_{n+1})$, 
where $c$ is a positive integer. It is well-known fact that $V$
can be realized as a submodule of a tensor product of $c$ level one modules,  $L(\Lambda_{j_1})\otimes\ldots\otimes L(\Lambda_{j_c})$,
generated by a highest weight vector $v_{\Lambda}=v_{\Lambda_{j_1}}\otimes\ldots\otimes v_{\Lambda_{j_c}}$.
Consider  an action of $x_{\alpha_{i}}^{+}(z_{1})\ldots x_{\alpha_{i}}^{+}(z_{m})$ on V. A formula for such an action consists of the 
summands $x_{\alpha_{i}}^{+(l_1;c)}(z_1)\cdots x_{\alpha_{i}}^{+(l_m;c)}(z_m)$ (see (\ref{103_coproduct})) whose components are the products of the operators $x_{\alpha_{i}}^{+}(z_{r}q^{l-1})$ and $\phi_{i}(z_{s}q^{l-1/2})$.
Our goal is to find a polynomial in
$$\mathcal{P}_{m}:=\mathbb{C}(q)\left[\frac{z_{s}}{z_{r}} : r,s=1,2,\ldots,m, r<s\right]$$ 
such that a product of 
$x_{\alpha_{i}}^{+}(z_{1})\ldots x_{\alpha_{i}}^{+}(z_{m})$ and this polynomial is a 
$\mathcal{P}_{m}[z_1,\ldots,z_{m-1}]$-linear combination of summands whose each component is ordered as in (\ref{order}).

\begin{lem}\label{type1}
For every  positive integer $m$ we have
$$\left(\prod_{r=1}^{m-1}\prod_{s=r+1}^{m}
\left(1-q^{2}\frac{z_{s}}{z_{r}}\right)\right)x_{\alpha_{i}}^{+}(z_{1})\ldots x_{\alpha_{i}}^{+}(z_{m})\in \om(V,V((z_1,\ldots,z_m))).$$
\end{lem}

\begin{dok}
Recall the formulas   (\ref{f:3}) and (\ref{f:4}). 
They imply that by multiplying the operator
$x_{\alpha_{i}}^{+}(z_{1})\ldots x_{\alpha_{i}}^{+}(z_{m})$ by polynomials $1-q^{2}z_{s}/z_{r}$, $r<s$, 
we  obtain the order (\ref{order}) 
on each component of each summand $x_{\alpha_{i}}^{+(l_1;c)}(z_1)\cdots x_{\alpha_{i}}^{+(l_m;c)}(z_m)$. 
Hence an $l$-th component ($l=1,2,\ldots,c$) of each summand of
$$\left(\prod_{r=1}^{m-1}\prod_{s=r+1}^{m}
\left(1-q^{2}\frac{z_{s}}{z_{r}}\right)\right)x_{\alpha_{i}}^{+}(z_{1})\ldots x_{\alpha_{i}}^{+}(z_{m})$$
is an element of  $\om(L(\Lambda_{j_l}),L(\Lambda_{j_l})((z_1,\ldots,z_m)))$. We conclude that 
$$\left(\prod_{r=1}^{m-1}\prod_{s=r+1}^{m}
\left(1-q^{2}\frac{z_{s}}{z_{r}}\right)\right)x_{\alpha_{i}}^{+(l_1;c)}(z_1)\cdots x_{\alpha_{i}}^{+(l_m;c)}(z_m)\in \om(V,V((z_1,\ldots,z_m)))$$
thus proving the lemma.
\end{dok}

The statement of Lemma \ref{type1} and the discussion preceding the lemma allow us to  define an  operator $x^{+}_{m\alpha_{i}}(z)$ on $V$,
\begin{equation}\label{type_1_op.}
x^{+}_{m\alpha_{i}}(z):=\lim_{z_{p}\to zq^{2(p-1)}}\left(\prod_{r=1}^{m-1}\prod_{s=r+1}^{m}
\left(1-q^{2}\frac{z_{s}}{z_{r}}\right)\right)x_{\alpha_{i}}^{+}(z_{1})\ldots x_{\alpha_{i}}^{+}(z_{m})\in\om(V,V((z)))
\end{equation}
for every positive integer $m$.
Of course, for $m=1$  we have
$x^{+}_{1\alpha_{i}}(z)=x^{+}_{\alpha_{i}}(z)$.
Notice that
$$\lim_{z_{p}\to zq^{2(p-1)}}\prod_{r=1}^{m-1}\prod_{s=r+1}^{m}
\left(1-q^{2}\frac{z_{s}}{z_{r}}\right)=\prod_{r=1}^{m-1}\prod_{s=r+1}^{m}
\left(1-q^{2(1+s-r)}\right)\neq 0.$$

\begin{defn}\label{type_1_def}
For an integer $r$ and a positive integer $m$ define  
$$x^{+}_{m\alpha_{i}}(r):=\rez_{z}z^{m+r-1}x_{m\alpha_{i}}^{+}(z)\in\ndo(V).$$
We call $x^{+}_{m\alpha_{i}}(r)$ a quasi-particle of type $1$, color $i$, charge $m$ and degree $r$. 
\end{defn}
We have
$$x_{m\alpha_{i}}^{+}(z)=\sum_{r\in\mathbb{Z}}x^{+}_{m\alpha_{i}}(r)z^{-r-m}.$$


\subsection{Quasi-particles of type 2}\label{sub:order}
In this subsection we introduce the notion of a quasi-particle of type $2$ (of certain color, 
charge and degree) for every integrable highest weight module $V:=L(\Lambda)$ of $U_{q}(\widehat{\mathfrak{sl}}_{n+1})$. 
These operators will be constructed using Ding-Feigin operators (\ref{105_ding_feigin}) from \cite{DF}.

Fix $i=1,2,\ldots,n$. For every positive integer $m$  define an operator $\bar{x}_{m\alpha_{i}}^{+}(z)$ on $V$:
\begin{equation}\label{type2}
\bar{x}_{m\alpha_{i}}^{+}(z):=\bar{x}^{+}_{\alpha_{i}}(z)\bar{x}^{+}_{\alpha_{i}}(zq^{2})\cdots\bar{x}_{\alpha_{i}}^{+}(zq^{2(m-1)}).
\end{equation}

\begin{lem}\label{restricted}
$\bar{x}_{m\alpha_{i}}^{+}(z)\in\om(V,V((z))).$
\end{lem}

\begin{dok}
Since the operator $\bar{x}^{+}_{\alpha_{i}}(z)\in\om(V,V((z)))$ commutes with itself (Theorem \ref{106_ding_feigin}), 
expression (\ref{type2}) is obviously well-defined on $V$ and $\bar{x}_{m\alpha_{i}}^{+}(z)\in\om(V,V((z)))$. 
\end{dok}

\begin{defn}
For an integer $r$ and a positive integer $m$ define 
$$\bar{x}^{+}_{m\alpha_{i}}(r):=\rez_{z}z^{m+r-1}\bar{x}_{m\alpha_{i}}^{+}(z)\in\om(V).$$
We call $\bar{x}^{+}_{m\alpha_{i}}(r)$ a quasi-particle of type $2$, color $i$, charge $m$ and degree $r$. 
\end{defn}

The next statement is an easy consuequence of Theorem \ref{106_ding_feigin} and definition of $\bar{x}_{m\alpha_{i}}^{+}(z)$.

\begin{kor}\label{107_t2.comm.}
Let $i,j=1,2,\ldots,n$ such that $a_{ij}=0$. We have on $V$
\begin{align*}
&\bar{x}_{m\alpha_i}^{+}(z_1)\bar{x}_{k\alpha_i}^{+}(z_2)=\bar{x}_{k\alpha_i}^{+}(z_2)\bar{x}_{m\alpha_i}^{+}(z_1),\\
&\bar{x}_{m\alpha_i}^{+}(z_1)\bar{x}_{k\alpha_j}^{+}(z_2)=\bar{x}_{k\alpha_j}^{+}(z_2)\bar{x}_{m\alpha_i}^{+}(z_1).
\end{align*}
\end{kor} 

Denote by $\bar{\mathfrak{S}}_{W(\Lambda)}$ a set of monomials of type $2$ quasi-particles   that satisfy following  assumptions:
\begin{enumerate}
  \item Product of quasi-particles has its quasi-particle colors nonincreasing from right to left;
  \item Product of quasi-particles of the same color has its quasi-particle charges nonincreasing from right to left;
  \item Product of quasi-particles of the same color and charge has its quasi-particle degrees nonincreasing from right to left.
\end{enumerate}
The following lemma is an easy consequence of Lemma \ref{colors_2} and Corollary \ref{107_t2.comm.}.

\begin{lem}\label{simple_lemma}
For a given dominant integral highest weight $\Lambda$ the set 
$$\left\{bv_{\Lambda}\hspace{2pt}\left|\hspace{2pt}\right. b\in\bar{\mathfrak{S}}_{W(\Lambda)}\right\}$$ 
spans  the principal subspace $W(\Lambda)$. 
\end{lem}

 Georgiev introduced in \cite{G} a strict linear (lexicographic) order ``$<$'' and a strict partial order ``$\prec$'' on a certain 
set of monomials of quasi-particles. We will recall his terminology and then we will  apply the above-mentioned orders  on  certain subsets of $\bar{\mathfrak{S}}_{W(\Lambda)}$.

For given $r_n,\ldots,r_1\in\mathbb{Z}_{\geq 0}$,  $r:=\sum_{s=1}^{n}r_{s}$, 
consider color-ordered sequences of $r_{n}$ integers of color $n$, \ldots, $r_2$ integers of color $2$, $r_1$ integers of color $1$:
$$m_{r}\leq \ldots\leq m_{\sum_{s=1}^{n-1}r_{s}+1},\quad m_{\sum_{s=1}^{n-1}r_{s}}\leq \quad\ldots\quad\leq m_{r_{1}+1},\quad m_{r_1}\leq \ldots\leq m_1,$$
such that only the entries of the same color are nonincreasing from right to left.
For two such sequences we write
\begin{equation}\label{ordering}
(m_{r},\ldots,m_{1})<(m^{'}_{r},\ldots,m^{'}_{1})
\end{equation}
if there exists $l\in\mathbb{Z}$, $1\leq l\leq r$, such that
$$m_{1}=m^{'}_{1},\quad m_{2}=m^{'}_{2},\quad\ldots.\quad m_{l-1}=m^{'}_{l-1}\quad\textrm{and}\quad m_{l}<m^{'}_{l}.$$
We write
\begin{equation}\label{ordering2}
(m_{r},\ldots,m_{1})\prec (m^{'}_{r},\ldots,m^{'}_{1})
\end{equation}
if there exists $l\in\mathbb{Z}$, $1\leq l\leq r$, such that
\begin{align*}
&m_{l}+m_{l-1}+\ldots +m_{2}+m_{1}< m^{'}_{l}+m^{'}_{l-1}+\ldots +m^{'}_{2}+m^{'}_{1},\\
&m_{k}+m_{k-1}+\ldots +m_{2}+m_{1}\leq m^{'}_{k}+m^{'}_{k-1}+\ldots +m^{'}_{2}+m^{'}_{1}\quad\textrm{for } k=1,2,\ldots,r.
\end{align*}

Fix a color $i=1,2,\ldots,n$.  We define a \emph{charge-type} of a monomial
$$b_{i}:=\bar{x}^{+}_{m_{r^{(1)}}\alpha_{i}}(s_{r^{(1)}})\cdots \bar{x}^{+}_{m_{1}\alpha_{i}}(s_{1})\in \bar{\mathfrak{S}}_{W(\Lambda)},$$
(consisting of quasi-particles of color $i$) as an $r^{(1)}$-tuple
$$(m_{r^{(1)}},m_{r^{(1)}-1},\ldots,m_{1}).$$ 
We define a \emph{dual-charge-type} of a monomial $b_{i}$ as an $m_1$-tuple
$$(r^{(1)},r^{(2)},\ldots,r^{(m_1)})$$
if it is built out of $r^{(1)}-r^{(2)}$ quasi-particles of charge $1$, $r^{(2)}-r^{(3)}$ quasi-particles of charge $2$, \ldots, $r^{(m_1)}$
quasi-particles of charge $m_1$.

Fix a monomial $b\in \bar{\mathfrak{S}}_{W(\Lambda)}$,
$b=b_{n}\cdots b_{2}b_{1}$,
where $b_{i}\in \bar{\mathfrak{S}}_{W(\Lambda)}$ is a monomial consisting of quasi-particles of color $i$.
We define a \emph{color-charge-type} of $b$ as an $r$-tuple
\begin{equation*}
(m_{r_{n}^{(1)},n},\ldots,m_{1,n};\ldots;m_{r_{1}^{(1)},1},\ldots,m_{1,1}),
\end{equation*}
if $b_{i}$ is of a charge-type  $(m_{r_{i}^{(1)},i},\ldots,m_{1,i})$ for every $i=1,2,\ldots,n$,
where
$$0<m_{r^{(1)}_{i},i}\leq\ldots\leq m_{2,i}\leq m_{1,i}\quad\textrm{for } i=1,2,\ldots,n,$$
and $r=r_{1}^{(1)}+r_{2}^{(1)}+\ldots+r_{n}^{(1)}$. 
We define a \emph{color-type} of $b$ as an $n$-tuple
$$(m_{n},\ldots,m_{1}),$$
where $m_{i}:=\sum_{s=1}^{r_{i}^{(1)}}m_{s,i}$ for $i=1,2,\ldots,n$.
We define a \emph{color-dual-charge-type} of $b$ as an $m$-tuple
\begin{equation*}
(r_{n}^{(1)},\ldots,r_{n}^{(m_{1,n})};\ldots;r_{1}^{(1)},\ldots,r_{1}^{(m_{1,1})}),
\end{equation*}
if $b_{i}$ is of a dual-charge-type $(r_{i}^{(1)},\ldots,r_{i}^{(m_{1,i})})$ for every $i=1,2,\ldots,n$,
where
$$r_{i}^{(1)}\geq r_{i}^{(2)}\geq\ldots\geq r_{i}^{(m_{1,i})}\geq 0\quad\textrm{for } i=1,2,\ldots,n$$
and $m=m_{1,1}+m_{1,2}+\ldots+m_{1,n}$. We will also say that the corresponding operator
$$\bar{x}^{+}_{m_{r_{n}^{(1)},n}\alpha_{n}}(z_{{r_{n}^{(1)}},n})\cdots \bar{x}^{+}_{m_{1,1}\alpha_{1}}(z_{1,1})$$
has the above color-charge-type, color-type and color-dual-charge-type.
We define a \emph{color-degree-type} of $b$ as an $n$-tuple $(l_{n},\ldots,l_2,l_{1})$,
where $l_{i}$, $i=1,2,\ldots,n$, is the sum of degrees of 
all the quasi-particles of color $i$ in $b$.

Fix a color-type  $(m_{n},\ldots,m_1)$. Now we can define a linear order ``$<$'' and a partial order ``$\prec$'' 
on the subset of $\bar{\mathfrak{S}}_{W(\Lambda)}$ that consists of monomials
of color-type $(m_{n},\ldots,m_1)$. The order ``$<$'' is defined as follows: 
First apply Definition (\ref{ordering}) to the color-charge-types of the two monomials $b$ and $b^{'}$;
if the color-charge-types are the same, aply (\ref{ordering}) to the degree sequences of the two monomials.
The partial order ``$\prec$'' is defined as follows: We write $b\prec b^{'}$ if $b<b^{'}$ and in addition
$$(l_{n},\ldots,l_{2},l_{1})\prec (l_{n}^{'},\ldots,l_{2}^{'},l_{1}^{'}),$$
(see (\ref{ordering2})),  where $(l_{n},l_{n-1},\ldots,l_{1})$ and $(l_{n}^{'},\ldots,l_{2}^{'},l_{1}^{'})$
are color-degree-types of quasi-particles $b$ and $b^{'}$ respectively.


\subsection{The main theorem} We shall consider only the principal subspaces $W(\Lambda)$ associated with the highest weights  $\Lambda\in \hat{P}^{+}$ of the form 
\begin{equation}\label{weight}
\Lambda=c_{0}\Lambda_{0}+c_{j}\Lambda_{j},
\end{equation}
such that $c_{0},c_{j}\in\mathbb{Z}_{\geq 0}$, $c_{0}+c_j>0$, $j=1,2,\ldots,n$.
The level $c$ of the weight $\Lambda$  equals $c_{0}+c_{j}$.
For every $j=1,2,\ldots,n$ define\begin{equation}\label{indices}j_{s}:=\begin{cases}0, & \textrm{if }s=1,2,\ldots ,c_{0},\\j, & \textrm{if }s=c_{0}+1,\ldots ,c_{0}+c_{j}.\end{cases}\end{equation}

For a given highest weight $\Lambda\in \hat{P}^{+}$ define a
set $\mathfrak{B}_{W(\Lambda)}$ of monomials of quasi-particles of type $1$,
\begin{align}\label{t2_basis}
\mathfrak{B}_{W(\Lambda)}&:=\bigcup_
{  \substack{   0\leq m_{r_{n}^{(1)}.n}\leq\ldots\leq m_{1,n}\leq c\\ \cdots\\0\leq m_{r_{1}^{(1)},1}\leq\ldots\leq m_{1,1}\leq c    }   }  \\
&\left\{x^{+}_{m_{r_{n}^{(1)},n}\alpha_{n}}(l_{r_{n}^{(1)},n})\cdots x^{+}_{m_{1,n}\alpha_{n}}(l_{1,n})\cdots x^{+}_{m_{r_{1}^{(1)},1}\alpha_{1}}(l_{r_{1}^{(1)},1})\cdots x^{+}_{m_{1,1}\alpha_{1}}(l_{1,1})\hspace{5pt}\Big|\right.\nonumber\\
&\qquad\Big|\hspace{5pt} l_{r,i}\leq \sum_{s=1}^{r_{i-1}^{(1)}}\min\left\{m_{r,i},m_{s,i-1}\right\}-\sum_{s=1}^{m_{r,i}}\delta_{ij_{s}}-\sum_{m_{t,i}>m_{r,i}}2m_{r,i}-m_{r,i},\nonumber\\
&\qquad \hspace{10pt}l_{r+1,i}\leq l_{r,i}-2m_{r,i}\textrm{ if }m_{r+1,i}=m_{r,i},\nonumber\\
&\qquad\bigg.\hspace{10pt}\textrm{for all } l_{r,i}\in\mathbb{Z},\hspace{3pt} i=1,2,\ldots,n,\hspace{3pt} r=1,2,\ldots,r_{i}^{(1)}\bigg\},\nonumber
\end{align}
and a set $\bar{\mathfrak{B}}_{W(\Lambda)}$ of monomials of quasi-particles of type $2$,
\begin{align}\label{t1_basis}
\bar{\mathfrak{B}}_{W(\Lambda)}&:=\bigcup_
{  \substack{   0\leq m_{r_{n}^{(1)},n}\leq\ldots\leq m_{1,n}\leq c\\ \cdots\\0\leq m_{r_{1}^{(1)},1}\leq\ldots\leq m_{1,1}\leq c    }   }  \\
&\left\{\bar{x}^{+}_{m_{r_{n}^{(1)},n}\alpha_{n}}(l_{r_{n}^{(1)},n})\cdots\bar{x}^{+}_{m_{1,n}\alpha_{n}}(l_{1,n})\cdots\bar{x}^{+}_{m_{r_{1}^{(1)},1}\alpha_{1}}(l_{r_{1}^{(1)},1})\cdots\bar{x}^{+}_{m_{1,1}\alpha_{1}}(l_{1,1})\hspace{5pt}\Big|\right.\nonumber\\
&\qquad\Big|\hspace{5pt} l_{r,i}\leq \sum_{s=1}^{r_{i-1}^{(1)}}\min\left\{m_{r,i},m_{s,i-1}\right\}-\sum_{s=1}^{m_{r,i}}\delta_{ij_{s}}-\sum_{m_{t,i}>m_{r,i}}2m_{r,i}-m_{r,i},\nonumber\\
&\qquad \hspace{10pt}l_{r+1,i}\leq l_{r,i}-2m_{r,i}\textrm{ if }m_{r+1,i}=m_{r,i},\nonumber\\
&\qquad\bigg.\hspace{10pt}\textrm{for all } l_{r,i}\in\mathbb{Z},\hspace{3pt} i=1,2,\ldots,n,\hspace{3pt} r=1,2,\ldots,r_{i}^{(1)}\bigg\}.\nonumber
\end{align}

The following theorem is the main result of this paper. Its proof will be given in the next two sections.

\noindent\textbf{Theorem \ref{main}}\textsl{ For a given highest weight $\Lambda=c_{0}\Lambda_{0}+c_{j}\Lambda_{j}$ as in (\ref{weight}) the sets 
$$\left\{bv_{\Lambda}\hspace{2pt}\left|\hspace{2pt}\right. b\in\mathfrak{B}_{W(\Lambda)}\right\}\qquad\textrm{and} \qquad \left\{bv_{\Lambda}\hspace{2pt}\left|\hspace{2pt}\right. b\in\bar{\mathfrak{B}}_{W(\Lambda)}\right\}$$ 
form the bases for the principal subspace $W(\Lambda)$. }


\section{The spanning set of \texorpdfstring{$W(\Lambda)$}{W(Lambda)}}


\subsection{Quantum integrability}
Fix $i=1,2,\ldots,n$.
In \cite{DM} J. Ding and T. Miwa discovered a condition of quantum integrability for level $c$ integrable $U_{q}(\widehat{\mathfrak{sl}}_{n+1})$-modules:
\begin{equation}\label{ding_miwa}
x_{\alpha_{i}}^{+}(z_{1})x_{\alpha_{i}}^{+}(z_{2})\cdots x_{\alpha_{i}}^{+}(z_{c+1})=0\quad\textrm{if}\quad z_{1}/z_{2}=z_{2}/z_{3}=\ldots=z_{c}/z_{c+1}=q^{-2}.
\end{equation}
By using the above result  J. Ding and B. Feigin (\cite{DF}) proved that on every level $c$ integrable $U_{q}(\widehat{\mathfrak{sl}}_{n+1})$-module
\begin{equation}\label{ding_feigin_3}
\bar{x}_{\alpha_{i}}^{+}(z_{1})\bar{x}_{\alpha_{i}}^{+}(z_{2})\cdots \bar{x}_{\alpha_{i}}^{+}(z_{c+1})=0\quad\textrm{if}\quad z_{1}/z_{2}=z_{2}/z_{3}=\ldots=z_{c}/z_{c+1}=q^{-2}.
\end{equation}
In the following proposition we formulate the relations of quantum integrability in terms of the operators $x_{m\alpha_{i}}^{+}(z)$ and $\bar{x}_{m\alpha_{i}}^{+}(z)$.
\begin{pro}\label{q_integrability}
 On every level $c$ integrable highest weight module 
$$x_{(c+1)\alpha_{i}}^{+}(z)=0\qquad\textrm{and}\qquad \bar{x}_{(c+1)\alpha_{i}}^{+}(z)=0.$$
\end{pro}

\begin{dok}
The first equality, $x_{(c+1)\alpha_{i}}^{+}(z)=0$, can be proved in the same way as (\ref{ding_miwa}) in \cite{DM} 
while the second equality, $\bar{x}_{(c+1)\alpha_{i}}^{+}(z)=0$, is  equivalent to (\ref{ding_feigin_3}).
\end{dok}

Notice that the sets $\mathfrak{B}_{W(\Lambda)}$ and $\bar{\mathfrak{B}}_{W(\Lambda)}$ 
contain only monomials of quasi-particles of charge less than or equal to $c$. This is the direct consequence of the above
proposition. In fact, we can formulate a simple refinement of Lemma \ref{simple_lemma}.
Denote by $\bar{\mathfrak{S}}_{W(\Lambda)}^{(c)}$ a subset of $\bar{\mathfrak{S}}_{W(\Lambda)}$ containing all  the monomials 
in $\bar{\mathfrak{S}}_{W(\Lambda)}$ that are built out of  quasi-particles of charge less than or equal to $c$.

\begin{lem}\label{simple_lemma2}
For a given level $c$ dominant integral  highest weight $\Lambda$  the set 
$$\left\{bv_{\Lambda}\hspace{2pt}\left|\hspace{2pt}\right. b\in\bar{\mathfrak{S}}_{W(\Lambda)}^{(c)}\right\}$$ 
spans  the principal subspace $W(\Lambda)$. 
\end{lem}


\subsection{Relations among quasi-particles of adjacent colors}
Let $i=1,2,\ldots,n$ and $m\in\mathbb{Z}_{>0}$. 
By applying the formula (\ref{103_coproduct}) on the operator $x_{m\alpha_{i}}^{+}(z)$ (defined by (\ref{type_1_op.})) we get
\begin{align*}
\Delta^{(c-1)}(x_{m\alpha_{i}}^{+}(z))=\hspace{-7pt}
\sum_{l_1,\ldots,l_m =1}^{c}\hspace{-4pt}\left(\hspace{-1pt}\lim_{z_{p}\to zq^{2(p-1)}}\left(\prod_{r=1}^{m-1}\prod_{s=r+1}^{m}
\left(1-q^{2}\frac{z_{s}}{z_{r}}\right)\hspace{-3pt}\right) x_{\alpha_{i}}^{+(l_1;c)}(z_{1})\ldots x_{\alpha_{i}}^{+(l_m;c)}(z_{m})\right)\hspace{-2pt}.
\end{align*}
Set
$$x_{m\alpha_{i}}^{+(l_1,\ldots,l_m;c)}(z):=\lim_{z_{p}\to zq^{2(p-1)}}\left(\prod_{r=1}^{m-1}\prod_{s=r+1}^{m}
\left(1-q^{2}\frac{z_{s}}{z_{r}}\right)\right) x_{\alpha_{i}}^{+(l_1;c)}(z_{1})\ldots x_{\alpha_{i}}^{+(l_m;c)}(z_{m}).$$
By using  (\ref{f:1}) we can analogously introduce the preceding notation for type $2$ quasi-particles.
We have
\begin{align}\label{cons:calc}
\bar{x}_{m\alpha_{i}}^{+}(z)&=\bar{x}_{\alpha_{i}}^{+}(z)\bar{x}_{\alpha_{i}}^{+}(zq^{2})\cdots\bar{x}_{\alpha_{i}}^{+}(zq^{2(m-1)})
=\lim_{z_{p}\to zq^{2(p-1)}}\bar{x}_{\alpha_{i}}^{+}(z_{1})\bar{x}_{\alpha_{i}}^{+}(z_{2})\cdots\bar{x}_{\alpha_{i}}^{+}(z_{m})\\
&=\lim_{z_{p}\to zq^{2(p-1)}}x_{\alpha_{i}}^{+}(z_{1})k_{i}^{+}(z_{1})x_{\alpha_{i}}^{+}(z_{2})k_{i}^{+}(z_{2})\cdots x_{\alpha_{i}}^{+}(z_{m})k_{i}^{+}(z_{m})\nonumber\\
&=\lim_{z_{p}\to zq^{2(p-1)}}
\left(\prod_{r=1}^{m-1}\prod_{s=r+1}^{m}\left(\frac{1-q^{2}\frac{z_{s}}{z_{r}}}{1-\frac{z_{s}}{z_{r}}}\right)\right)
x_{\alpha_{i}}^{+}(z_{1})\cdots x_{\alpha_{i}}^{+}(z_{m})k_{i}^{+}(z_{1})\cdots k_{i}^{+}(z_{m}).\nonumber
\end{align}
Set
\begin{align*}
&\bar{x}_{m\alpha_{i}}^{+(l_1,\ldots,l_m;c)}(z):=\\
&\hspace{20pt}\lim_{z_{p}\to zq^{2(p-1)}}\left(\prod_{r=1}^{m-1}\prod_{s=r+1}^{m}\left(\frac{1-q^{2}\frac{z_{s}}{z_{r}}}{1-\frac{z_{s}}{z_{r}}}\right)\right)
 x_{\alpha_{i}}^{+(l_1;c)}(z_{1})\ldots x_{\alpha_{i}}^{+(l_m;c)}(z_{m})k_{i}^{+}(z_{1})\cdots k_{i}^{+}(z_{m}).
\end{align*}

\begin{lem}\label{x_order}
On every level $c$ integrable highest weight module we have
\begin{enumerate}
\item If there exist integers $r,s=1,2,\ldots,c$, $r\neq s$, such that $l_{r}=l_{s}$,  then
 $$x_{m\alpha_{i}}^{+(l_1,\ldots,l_m;c)}(z)=0\qquad\textrm{and}\qquad \bar{x}_{m\alpha_{i}}^{+(l_1,\ldots,l_m;c)}(z)=0;$$
\item If there exists an integer $r=1,2,\ldots,c-1$ such that $l_{r}<l_{r+1}$, then
 $$x_{m\alpha_{i}}^{+(l_1,\ldots,l_m;c)}(z)=0\qquad\textrm{and}\qquad \bar{x}_{m\alpha_{i}}^{+(l_1,\ldots,l_m;c)}(z)=0.$$
\end{enumerate}
\end{lem}

\begin{dok}
(1) Fix $m$ integers $ l_{1},\ldots, l_{m}=1,2,\ldots,c$ and assume that there exist $r,s=1,2,\ldots,m$ such that $l_{r}=l_{s}$ and $ r< s$.

If $r+1=s$ the equality $x_{m\alpha_{i}}^{+(l_1,\ldots,l_m;c)}(z)=0$ follows from quantum integrability (see Proposition \ref{q_integrability}).
If $r+1<s$ notice that there exists an integer $u$ such that $r\leq r+u< s$ and $l_{r+u}\leq l_{r+u+1}$ (on the contrary we would
have $l_{r}>l_{r+1}>\ldots >l_{s-1}>l_s$ and $l_r =l_s$). 
If $l_{r+u}= l_{r+u+1}$ we can conlude that $x_{m\alpha_{i}}^{+(l_1,\ldots,l_m;c)}(z)=0$ in the same way as we did above. 
If $l_{r+u}< l_{r+u+1}$ we observe, inside of a limit $\lim_{z_{p}\to zq^{2(p-1)}}$, the $l_{r+u}$-th tensor component of $x_{m\alpha_{i}}^{+(l_1,\ldots,l_m;c)}(z)$:
$$\ldots x_{\alpha_{i}}^{+}(z_{r+u}q^{l_{r+u}-1})\phi_{i}(z_{r+u+1}q^{l_{r+u}-1/2})\ldots.$$
By applying   (\ref{f:4}) we get
$$\ldots q^{2}\left(1-q^{-2}\frac{z_{r+u+1}}{z_{r+u}}\right)\phi_{i}(z_{r+u+1}q^{l_{r+u}-1/2})x_{\alpha_{i}}^{+}(z_{r+u}q^{l_{r+u}-1})\ldots.$$
The limit $\lim_{z_{p}\to zq^{2(p-1)}}$ annihilates the $l_{r+u}$-th tensor component because of the obtained factor $(1-q^{-2}z_{r+u+1}/z_{r+u})$ 
so we conclude that $x_{m\alpha_{i}}^{+(l_1,\ldots,l_m;c)}(z)=0$. It is important to emphasize that the formulas 
(\ref{f:3}) and (\ref{f:4}) guarantee that in  
$x_{\alpha_{i}}^{+(l_1;c)}(z_{1})\ldots x_{\alpha_{i}}^{+(l_m;c)}(z_{m})$
will not appear any term that would cancel the annihilating term $(1-q^{-2}z_{r+u+1}/z_{r+u})$.

Now we can prove that under the above assumptions $\bar{x}_{m\alpha_{i}}^{+(l_1,\ldots,l_m;c)}(z)=0$. Define a Laurent polynomial
$p(z_1,\ldots,z_m):=\prod_{r=1}^{m-1}\prod_{s=r+1}^{m}(1-z_{s}/z_{r})$. The equality 
$x_{m\alpha_{i}}^{+(l_1,\ldots,l_m;c)}(z)=0$ implies 
\begin{align*}
\lim_{z_{p}\to zq^{2(p-1)}}\Bigg(&p(z_1,\ldots,z_m)\left(\prod_{r=1}^{m-1}\prod_{s=r+1}^{m}\left(\frac{1-q^{2}\frac{z_{s}}{z_{r}}}{1-\frac{z_{s}}{z_{r}}}\right)\right)\\
&\cdot x_{\alpha_{i}}^{+(l_1;c)}(z_{1})\ldots x_{\alpha_{i}}^{+(l_m;c)}(z_{m})k_{i}^{+}(z_{1})\cdots k_{i}^{+}(z_{m})\Bigg)=0.
\end{align*}
 Since $\lim_{z_{p}\to zq^{2(p-1)}}p(z_1,\ldots,z_m)\neq 0$ we conclude that
\begin{align*}
&\bar{x}_{m\alpha_{i}}^{+(l_1,\ldots,l_m;c)}(z)\\
&=\lim_{z_{p}\to zq^{2(p-1)}}\left(\prod_{r=1}^{m-1}\prod_{s=r+1}^{m}\left(\frac{1-q^{2}\frac{z_{s}}{z_{r}}}{1-\frac{z_{s}}{z_{r}}}\right)\right)
 x_{\alpha_{i}}^{+(l_1;c)}(z_{1})\ldots x_{\alpha_{i}}^{+(l_m;c)}(z_{m})k_{i}^{+}(z_{1})\cdots k_{i}^{+}(z_{m})=0.
 \end{align*} 
 
 (2) The second statement of the lema is an easy consequence of the proof of the first one. The assumption
 that there exists an integer $r=1,2,\ldots,c-1$ such that $l_{r}<l_{r+1}$ implies that, 
 inside of a limit $\lim_{z_{p}\to zq^{2(p-1)}}$,  on the $l_{r}$-th component  of $x_{m\alpha_{i}}^{+(l_1,\ldots,l_m;c)}(z)$ we have
 $$\ldots x_{\alpha_{i}}^{+}(z_{r}q^{l_{r}-1})\phi_{i}(z_{r+1}q^{l_{r}-1/2})\ldots.$$
 Now we can proceed as in the first part of the proof.
\end{dok}

By putting $m=c$ in the above lemma we are getting the quasi-particle analogue of Lemma 2.7. in \cite{DF2}. 
The proof technique we used is quite similar to the one in that paper.

Fix an integer $i=2,3,\ldots,n$.

\begin{lem}\label{tech1}
On every level $1$ integrable highest weight module we have
 \begin{align*} 
 &x_{\alpha_{i}}^{+}(z_{1}q^{2(r-1)})\phi_{i-1}(z_{2}q^{\frac{1}{2}})\phi_{i-1}(z_{2}q^{\frac{5}{2}})\ldots\phi_{i-1}(z_{2}q^{2(s-2)+\frac{1}{2}}) x_{\alpha_{i-1}}^{+}(z_{2}q^{2(s-1)})\\
 &\hspace{20pt}=q^{1-s}\frac{1}{z_{1}}\frac{1}{1-q^{1-2r}\frac{z_{2}}{z_{1}}}\phi_{i-1}(z_{2}q^{\frac{1}{2}})\phi_{i-1}(z_{2}q^{\frac{5}{2}})\ldots\phi_{i-1}(z_{2}q^{2(s-2)+\frac{1}{2}})\\
 &\hspace{40pt}\cdot:x_{\alpha_i}^{+}(z_{1}q^{2(r-1)})x_{\alpha_{i-1}}^{+}(z_{2}q^{2(s-1)}):
 \end{align*}
for $r,s\in\mathbb{Z}_{>0}$.
\end{lem}

\begin{dok}
Formula (\ref{f:6}) implies
 \begin{align}\label{imp:eq}
 &x_{\alpha_{i}}^{+}(z_{1}q^{2(r-1)})\phi_{i-1}(z_{2}q^{\frac{1}{2}})\phi_{i-1}(z_{2}q^{\frac{5}{2}})\ldots\phi_{i-1}(z_{2}q^{2(s-2)+\frac{1}{2}})\\
 &\hspace{15pt}=q^{1-s}\left(\prod_{t=2}^{s}\frac{1-q^{2(t-r)-1}\frac{z_{2}}{z_{1}}}{1-q^{2(t-r)-3}\frac{z_{2}}{z_{1}}}\right)\phi_{i-1}(z_{2}q^{\frac{1}{2}})\phi_{i-1}(z_{2}q^{\frac{5}{2}})\ldots\phi_{i-1}(z_{2}q^{2(s-2)+\frac{1}{2}})x_{\alpha_i}^{+}(z_{1}q^{2(r-1)})\nonumber\\
 &\hspace{15pt}=q^{1-s}\frac{1-q^{2(s-r)-1}\frac{z_{2}}{z_{1}}}{1-q^{1-2r}\frac{z_{2}}{z_{1}}}\phi_{i-1}(z_{2}q^{\frac{1}{2}})\phi_{i-1}(z_{2}q^{\frac{5}{2}})\ldots\phi_{i-1}(z_{2}q^{2(s-2)+\frac{1}{2}})x_{\alpha_i}^{+}(z_{1}q^{2(r-1)}).\nonumber
 \end{align}
 Finally, by using (\ref{f:5}) we get
\begin{align*}
 &x_{\alpha_{i}}^{+}(z_{1}q^{2(r-1)})\phi_{i-1}(z_{2}q^{\frac{1}{2}})\phi_{i-1}(z_{2}q^{\frac{5}{2}})\ldots\phi_{i-1}(z_{2}q^{2(s-2)+\frac{1}{2}})x_{\alpha_{i-1}}^{+}(z_{2}q^{2(s-1)})\\
 &\hspace{20pt}=q^{1-s}\frac{1-q^{2(s-r)-1}\frac{z_{2}}{z_{1}}}{1-q^{1-2r}\frac{z_{2}}{z_{1}}}\phi_{i-1}(z_{2}q^{\frac{1}{2}})\phi_{i-1}(z_{2}q^{\frac{5}{2}})\ldots\phi_{i-1}(z_{2}q^{2(s-2)+\frac{1}{2}})\\
 &\hspace{45pt}\cdot x_{\alpha_i}^{+}(z_{1}q^{2(r-1)})x_{\alpha_{i-1}}^{+}(z_{2}q^{2(s-1)})\\
  &\hspace{20pt}=q^{1-s}\frac{1-q^{2(s-r)-1}\frac{z_{2}}{z_{1}}}{1-q^{1-2r}\frac{z_{2}}{z_{1}}}\phi_{i-1}(z_{2}q^{\frac{1}{2}})\phi_{i-1}(z_{2}q^{\frac{5}{2}})\ldots\phi_{i-1}(z_{2}q^{2(s-2)+\frac{1}{2}})\\
  &\hspace{45pt}\cdot\frac{1}{z_1}\frac{1}{1-q^{2(s-r)-1}\frac{z_2}{z_1}}:x_{\alpha_i}^{+}(z_{1}q^{2(r-1)})x_{\alpha_{i-1}}^{+}(z_{2}q^{2(s-1)}):\\
  &\hspace{20pt}=q^{1-s}\frac{1}{z_1}\frac{1}{1-q^{1-2r}\frac{z_{2}}{z_{1}}}\phi_{i-1}(z_{2}q^{\frac{1}{2}})\phi_{i-1}(z_{2}q^{\frac{5}{2}})\ldots\phi_{i-1}(z_{2}q^{2(s-2)+\frac{1}{2}})\\
  &\hspace{45pt}\cdot:x_{\alpha_i}^{+}(z_{1}q^{2(r-1)})x_{\alpha_{i-1}}^{+}(z_{2}q^{2(s-1)}):.
 \end{align*}
\end{dok}

The next lemma is a consequence of (\ref{imp:eq}).

\begin{lem}\label{tech2}
On every level $1$ integrable highest weight module we have
 \begin{align*} 
 &x_{\alpha_i}^{+}(z_{1}q^{2(r-1)})\phi_{i-1}(z_{2}q^{\frac{1}{2}})\phi_{i-1}(z_{2}q^{\frac{5}{2}})\ldots\phi_{i-1}(z_{2}q^{2(s-1)+\frac{1}{2}})\\
 &\hspace{30pt}=\frac{1-q^{2(s-r)+1}\frac{z_{2}}{z_{1}}}{1-q^{1-2r}\frac{z_{2}}{z_{1}}}\phi_{i-1}(z_{2}q^{\frac{1}{2}})\phi_{i-1}(z_{2}q^{\frac{5}{2}})\ldots\phi_{i-1}(z_{2}q^{2(s-2)+\frac{1}{2}})x_{\alpha_{i}}^{+}(z_{1}q^{2(r-1)})
 \end{align*}
 for $r,s\in\mathbb{Z}_{>0}$.
\end{lem}

Fix a positive integer $k$.
In order to simplify some of the formulas in this paragraph, we will omit operator variables in some of them. If $x_{m\alpha_{i}}^{+(l_1,\ldots,l_m;c)}(z)\neq 0$ then, by Lemma \ref{x_order}, every tensor component of 
$x_{m\alpha_{i}}^{+(l_1,\ldots,l_m;c)}(z)$ consists of zero or more   operators $\phi_{i}$ and of at 
most one operator $x_{\alpha_{i}}^{+}$  that is positioned on the right of all  the operators $\phi_{i}$.
Therefore,  a random tensor component consists of the operators positioned in the following order: 
$\phi_{i}\ldots \phi_{i}x_{\alpha_{i}}^{+}$; 
possibly without any operator $\phi_{i}$ or without an operator  $x_{\alpha_{i}}^{+}$ or without any of the 
operators, being equal to $1$.
If $x_{m\alpha_{i}}^{+(l_1,\ldots,l_m;c)}(z_1)x_{k\alpha_{i-1}}^{+(t_1,\ldots,t_k;c)}(z_2)\neq 0$ then  
every tensor component of $x_{m\alpha_{i}}^{+(l_1,\ldots,l_m;c)}(z_1)x_{k\alpha_{i-1}}^{+(t_1,\ldots,t_k;c)}(z_2)$ 
consists of the operators positioned in the following order:
$\phi_{i}\ldots \phi_{i}x_{\alpha_{i}}^{+}\phi_{i-1}\ldots \phi_{i-1}x_{\alpha_{i-1}}^{+}$;
possibly without some operators or without all  the operators,  being equal to $1$.
Notice that Lemma \ref{tech1} and Lemma  \ref{tech2} allow us to change the order of the operators on every tensor component. 
We will denote by $:x_{m\alpha_{i}}^{+(l_1,\ldots,l_m;c)}(z_1)x_{k\alpha_{i-1}}^{+(t_1,\ldots,t_k;c)}(z_2):$  
an operator whose every tensor component consists of the same operators as the corresponding tensor component
of $x_{m\alpha_{i}}^{+(l_1,\ldots,l_m;c)}(z_1)x_{k\alpha_{i-1}}^{+(t_1,\ldots,t_k;c)}(z_2)$ 
but with the operators on every tensor component positioned in the following order: 
\begin{equation}\label{op:order}
\phi_{i}\ldots \phi_{i}\phi_{i-1}\ldots \phi_{i-1}:x_{\alpha_{i}}^{+}x_{\alpha_{i-1}}^{+}:.
\end{equation}

\begin{lem}\label{tech3}
There exists a polynomial $B(z)\in\mathbb{C}(q^{1/2})[z]$ such that on every level $c$ integrable highest weight module
\begin{align*}
&z_{1}^{\min\left\{m,k\right\}}B(z_{2}/z_{1})x_{m\alpha_{i}}^{+(l_1,\ldots,l_m;c)}(z_1)x_{k\alpha_{i-1}}^{+(t_1,\ldots,t_k;c)}(z_2)\\
&\hspace{20pt}\in \mathbb{C}(q^{1/2})[z_1, z_2 /z_1]:x_{m\alpha_{i}}^{+(l_1,\ldots,l_m;c)}(z_1)x_{k\alpha_{i-1}}^{+(t_1,\ldots,t_k;c)}(z_2):\nonumber 
\end{align*}
for all integers $l_1,\ldots,l_m,t_1,\ldots,t_k=1,2,\ldots,c$.
The polynomial $B(z)$ is given by
\begin{equation}\label{polynomial}
B(z)=\begin{cases}
\displaystyle\prod_{r=1}^{m}\left(1-q^{1-2r}z\right), & \textrm{ if }m\leq k,\\
\displaystyle\prod_{r=m-k+1}^{m}\left(1-q^{1-2r}z\right), & \textrm{ if }m>k.
\end{cases}
\end{equation}
\end{lem}

\begin{dok}
Without loss of generality we can assume $l_1>l_2>\ldots>l_m$ and $t_1>t_2>\ldots>t_k$ because in the contrary Lemma \ref{x_order}
would imply $x_{m\alpha_{i}}^{+(l_1,\ldots,l_m;c)}(z_1)=0$ or $x_{k\alpha_{i-1}}^{+(t_1,\ldots,t_k;c)}(z_2)=0$.

Suppose $m\leq k$. For every $l=1,2,\ldots,c$, $l\neq l_{r}$, $r=1,2,\ldots,m$,  the $l$-th tensor component
of $x_{m\alpha_{i}}^{+(l_1,\ldots,l_m;c)}(z_1)x_{k\alpha_{i-1}}^{+(t_1,\ldots,t_k;c)}(z_2)$ consists of the operators
ordered as in (\ref{op:order}). On the remaining $m$ components we can apply Lemma \ref{tech1} or \ref{tech2} as appropriate. 
By multiplying the $l_r$-th component  by a factor $z_{1}(1-q^{1-2r}z_{2}/z_{1})$, where $r=1,2,\ldots,m$,
we order the operators of that component as in (\ref{op:order}). Thereby we proved
\begin{align}\label{case:1}
&z_{1}^{\min\left\{m,k\right\}}\prod_{r=1}^{m}\left(1-q^{1-2r}\frac{z_2}{z_1}\right)x_{m\alpha_{i}}^{+(l_1,\ldots,l_m;c)}(z_1)x_{k\alpha_{i-1}}^{+(t_1,\ldots,t_k;c)}(z_2)\\
&\hspace{20pt}\in \mathbb{C}(q^{1/2})[z_1, z_2 /z_1]:x_{m\alpha_{i}}^{+(l_1,\ldots,l_m;c)}(z_1)x_{k\alpha_{i-1}}^{+(t_1,\ldots,t_k;c)}(z_2):.\nonumber
\end{align}

Suppose $m> k$. Notice that by proceeding as in the case $m\leq k$ we can actually prove (\ref{case:1}). However, we want  $B(z)$ 
to be a polynomial of degree $\min\left\{m,k\right\}=k$ and (\ref{case:1}) gives us a polynomial of degree $m$. Assume that the 
statement of the lemma doesn't hold i.e. assume that 
$z_{1}^{\min\left\{m,k\right\}}B(z_{2}/z_{1})x_{m\alpha_{i}}^{+(l_1,\ldots,l_m;c)}(z_1)x_{k\alpha_{i-1}}^{+(t_1,\ldots,t_k;c)}(z_2)$
is not an element of
$\mathbb{C}(q^{1/2})[z_1, z_2 /z_1]:x_{m\alpha_{i}}^{+(l_1,\ldots,l_m;c)}(z_1)x_{k\alpha_{i-1}}^{+(t_1,\ldots,t_k;c)}(z_2):$
for some integers $l_1,\ldots,l_m,t_1,\ldots,t_k=1,2,\ldots,c$. 
This implies that when we order tensor components of
$z_{1}^{\min\left\{m,k\right\}}B(z_{2}/z_{1})x_{m\alpha_{i}}^{+(l_1,\ldots,l_m;c)}(z_1)x_{k\alpha_{i-1}}^{+(t_1,\ldots,t_k;c)}(z_2)$
as in (\ref{op:order}) (by using Lemma \ref{tech1} and Lemma \ref{tech2}), at least one factor 
$(1-q^{1-2r}z_{2}/z_{1})^{-1}$,
$r=1,2,\ldots,m-k$, appears so there exists an integer $s=1,2,\ldots,k$ such that $t_s\geq l_r$.

Suppose $t_k>l_{r+k}$. Denote the variables  in the following way
\begin{align*}
x_{m\alpha_i}^{+}(z_{1})&=\lim_{z^{'}_{p}\to z_{1}q^{2(p-1)}}\left(\prod_{r=1}^{m-1}\prod_{s=r+1}^{m}\left(1-q^{2}\frac{z_{s}^{'}}{z_{r}^{'}}\right)\right)x_{\alpha_i}^{+}(z^{'}_{1})x_{\alpha_i}^{+}(z^{'}_{2})\cdots x_{\alpha_i}^{+}(z^{'}_{m}),\\
x_{k\alpha_{i-1}}^{+}(z_{2})&=\lim_{z^{''}_{p}\to z_{2}q^{2(p-1)}}\left(\prod_{r=1}^{k-1}\prod_{s=r+1}^{k}\left(1-q^{2}\frac{z_{s}^{''}}{z_{r}^{''}}\right)\right)x_{\alpha_{i-1}}^{+}(z^{''}_{1})x_{\alpha_{i-1}}^{+}(z^{''}_{2})\cdots x_{\alpha_{i-1}}^{+}(z^{''}_{k}).
\end{align*}
Then, inside of the limits $\lim_{z^{'}_{p}\to z_{1}q^{2(p-1)}}$ and $\lim_{z^{''}_{p}\to z_{2}q^{2(p-1)}}$, 
we can carry out the following calculation on the $l_{r+k}$-th tensor component:
\begin{align*} 
 &x_{\alpha_i}^{+}(z_{r+k}^{'}q^{l_{1}-1})\phi_{i-1}(z_{1}^{''}q^{l_{1}-1/2})\ldots\phi_{i-1}(z_{k}^{''}q^{l_{1}-1/2})\\
 &\hspace{15pt}=q^{-k}\left(\prod_{s=1}^{k}\frac{1-q\frac{z_{s}^{''}}{z_{r+k}^{'}}}{1-q^{-1}\frac{z_{s}^{''}}{z_{r+k}^{'}}}\right)\phi_{i-1}(z_{1}^{''}q^{l_{1}-1/2})\ldots\phi_{i-1}(z_{k}^{''}q^{l_{1}-1/2})x_{\alpha_i}^{+}(z_{r+k}^{'}q^{l_{1}-1}).
 \end{align*}
 Since 
\begin{equation}\label{denom}
\lim_{\substack{ z^{'}_{p}\to z_{1}q^{2(p-1)}\\ z^{''}_{p}\to z_{2}q^{2(p-1)}}}\left(\prod_{s=1}^{k}\frac{1-q\frac{z_{s}^{''}}{z_{r+k}^{'}}}{1-q^{-1}\frac{z_{s}^{''}}{z_{r+k}^{'}}}\right)=\frac{1-q^{1-2r}\frac{z_{2}}{z_{1}}}{1-q^{1-2(r+k)}\frac{z_{2}}{z_{1}}}
\end{equation}
would have canceled the term $(1-q^{1-2r}z_{2}/z_{1})^{-1}$ and that would be in contradiction to our initial assumption, we conclude 
$t_k\leq l_{r+k}$. (Note that a denominator on the right hand side of (\ref{denom}) is a factor of $B(z_{2}/z_{1})$.)  Using the similar arguments and calculations we can prove $t_{k-1}\leq l_{r+k-1}$, $t_{k-2}\leq l_{r+k-2}$, \ldots
, $t_{2}\leq l_{r+2}$, $t_{1}\leq l_{r+1}$. Since $l_{r}>l_{r+1}$ we have $l_{r}>t_{s}$ for $s=1,2,\ldots,k$.  Contradiction!
\end{dok}

Lemma \ref{tech3} obviously implies:

\begin{kor}
On every level $c$ integrable highest weight module $L(\Lambda)$ we have
$$\left(\prod_{r=1}^{\min\left\{m,k\right\}}\left(z_1-q^{1-2(r+m-\min\left\{m,k\right\})}z_{2}\right)\right)x_{m\alpha_i}^{+}(z_{1})x_{k\alpha_{i-1}}^{+}(z_{2})\in\om(L(\Lambda),L(\Lambda)((z_{1},z_{2}))).$$
\end{kor}

As we will see later (see Definition \ref{quasi:comp}), the corollary 
actually states that an ordered pair $(x_{m\alpha_i}^{+}(z), x_{k\alpha_{i-1}}^{+}(z))$ is quasi compatible.


\begin{kor}
Let $B(z)$ be a polynomial defined by (\ref{polynomial}) and let $p$ be a positive integer, $p\leq \min\left\{m,k\right\}$.
Let $C(z)\in\mathbb{C}(q^{1/2})[z]$ be a polynomial such that
\begin{enumerate}
  \item $C(0)=1$,
  \item $B$ is divisible by $C$,
  \item For all integers $l_1,\ldots,l_m,t_1,\ldots,t_k=1,2,\ldots,c$
  \begin{align*}
&z_{1}^{p}C(z_{2}/z_{1})x_{m\alpha_{i}}^{+(l_1,\ldots,l_m;c)}(z_1)x_{k\alpha_{i-1}}^{+(t_1,\ldots,t_k;c)}(z_2)\\
&\hspace{20pt}\in \mathbb{C}(q^{1/2})[z_1, z_2 /z_1]:x_{m\alpha_{i}}^{+(l_1,\ldots,l_m;c)}(z_1)x_{k\alpha_{i-1}}^{+(t_1,\ldots,t_k;c)}(z_2):. 
\end{align*}
\end{enumerate}
Then $p= \min\left\{m,k\right\}$ and $C=B$.
\end{kor}

\begin{dok}
By carrying out a proof of Lemma \ref{tech3} for 
$(l_{1},l_2,\ldots,l_{m})=(m,m-1,\ldots,1)$ and $(t_1,t_2,\ldots,t_k)=(k,k-1,\ldots,1)$ we can se that
\begin{align*}
&z_{1}^{p}C(z_{2}/z_{1})x_{m\alpha_{i}}^{+(m,m-1,\ldots,2,1;c)}(z_1)x_{k\alpha_{i-1}}^{+(k,k-1,\ldots,2,1;c)}(z_2)\\
&\hspace{20pt}\in \mathbb{C}(q^{1/2})[z_1, z_2 /z_1]:x_{m\alpha_{i}}^{+(m,m-1,\ldots,2,1;c)}(z_1)x_{k\alpha_{i-1}}^{+(k,k-1,\ldots,2,1;c)}(z_2):
\end{align*}
implies $p= \min\left\{m,k\right\}$ and $C=B$.
\end{dok}

\begin{lem}\label{tech4} 
There exist a Taylor series $A(z)\in\mathbb{C}(q^{1/2})[[z]]$, $A(0)=1$, and a polynomial $B(z)\in\mathbb{C}(q^{1/2})[z]$ such that on every level $c$ integrable highest weight module
\begin{align*}
&z_{1}^{\min\left\{m,k\right\}}A(z_{2}/z_{1})B(z_{2}/z_{1})\bar{x}_{m\alpha_{i}}^{+(l_1,\ldots,l_m;c)}(z_1)\bar{x}_{k\alpha_{i-1}}^{+(t_1,\ldots,t_k;c)}(z_2)\\
&\hspace{20pt}\in \mathbb{C}(q^{1/2})[z_1, z_2 /z_1]:\bar{x}_{m\alpha_{i}}^{+(l_1,\ldots,l_m;c)}(z_1)\bar{x}_{k\alpha_{i-1}}^{+(t_1,\ldots,t_k;c)}(z_2):
\end{align*}
for all $l_1,\ldots,l_m,t_1,\ldots,t_k=1,2,\ldots,c$.
The polynomial $B(z)$ is given by (\ref{polynomial}). 
\end{lem}

\begin{dok}
Recall (\ref{cons:calc}). We can employ  (\ref{f:2}) in order  to shift the operators $k_{i}^{+}(z_{r}^{'})$, $r=1,2,\ldots,m$, in
\begin{align*}
&\bar{x}_{m\alpha_{i}}^{+}(z_1)\bar{x}_{k\alpha_{i-1}}^{+}(z_2)=\\
&\hspace{20pt}\lim_{\substack{z^{'}_{p}\to z_{1}q^{2(p-1)}\\z^{''}_{p}\to z_{2}q^{2(p-1)}}}\Bigg(
\left(\prod_{r=1}^{m-1}\prod_{s=r+1}^{m}\left(\frac{1-q^{2}\frac{z_{s}^{'}}{z_{r}^{'}}}{1-\frac{z_{s}^{'}}{z_{r}^{'}}}\right)\right)
\left(\prod_{r=1}^{k-1}\prod_{s=r+1}^{k}\left(\frac{1-q^{2}\frac{z_{s}^{''}}{z_{r}^{''}}}{1-\frac{z_{s}^{''}}{z_{r}^{''}}}\right)\right)\\
&\hspace{20pt}\cdot x_{\alpha_{i}}^{+}(z_{1}^{'})\cdots x_{\alpha_{i}}^{+}(z_{m}^{'})k_{i}^{+}(z_{1}^{'})\cdots k_{i}^{+}(z_{m}^{'})
x_{\alpha_{i-1}}^{+}(z_{1}^{''})\cdots x_{\alpha_{i-1}}^{+}(z_{k}^{''})k_{i-1}^{+}(z_{1}^{''})\cdots k_{i-1}^{+}(z_{k}^{''})\Bigg)
\end{align*}
to the right of the operators $x_{\alpha_{i-1}}^{+}(z_{s}^{''})$, $s=1,2,\ldots,k$. 
By doing this we will get a product of Taylor series in variables $z_{s}^{''}/z_{r}^{'}$, all of them having a constant term $1$.
A limit $z^{'}_{p}\to z_{1}q^{2(p-1)}$, $z^{''}_{p}\to z_{2}q^{2(p-1)}$ of this product is a Taylor series $D(z_2/z_1)$. 
Obviously $D(0)=1$ so we can set $A(z):=(D(z))^{-1}$. Now we can ensure, using Lemma \ref{tech3},  
that the series $A(z)$ would satisfy the statement of the lemma.
\end{dok}

In the next two results we will maintain the same as above notation: $B(z)$ will be  a polynomial given by (\ref{polynomial}) 
and $A(z)$ will be a Taylor series from 
Lemma \ref{tech4}.

\begin{lem}\label{tech5}
Let $\Lambda=c_{0}\Lambda_{0}+c_{j}\Lambda_{j}$ be a highest weight as in  (\ref{weight}) and let $v_{\Lambda}$ be a highest weight vector  
of $L(\Lambda)$.  We have
\begin{equation}\label{norm:exp}
z_{1}^{\min\left\{m,k\right\}}A(z_2/z_1)B(z_2/z_1)\bar{x}_{m\alpha_i}^{+}(z_{1})\bar{x}_{k\alpha_{i-1}}^{+}(z_{2})v_{\Lambda}
\in z_{1}^{\sum_{s=1}^{m}\delta_{ij_{s}}}z_{2}^{\sum_{s=1}^{k}\delta_{i-1\hspace{1pt}j_{s}}}W(\Lambda)[[z_{1},z_{2}]].
\end{equation}
\end{lem}

\begin{dok}
Since
$z_{1}^{\min\left\{m,k\right\}}A(z_{2}/z_{1})B(z_{2}/z_{1})\bar{x}_{m\alpha_{i}}^{+(l_1,\ldots,l_m;c)}(z_1)\bar{x}_{k\alpha_{i-1}}^{+(t_1,\ldots,t_k;c)}(z_2)$
is an element of a set
$\mathbb{C}(q^{1/2})[z_1, z_2 /z_1]:\bar{x}_{m\alpha_{i}}^{+(l_1,\ldots,l_m;c)}(z_1)\bar{x}_{k\alpha_{i-1}}^{+(t_1,\ldots,t_k;c)}(z_2):$
we conclude that, in (\ref{norm:exp}), on a vector $v_{\Lambda}=v_{\Lambda_{j_1}}\otimes\ldots\otimes v_{\Lambda_{j_c}}$ 
are first applied the operators $k_{i}^{+}$ and $k_{i-1}^{+}$. Of course, 
\begin{equation}\label{eqi:1}
k_{j}^{+}(z)v_{\Lambda}=v_{\Lambda}
\end{equation} 
for $j=1,2,\ldots,n$.
Next, a normal ordered product $:x_{\alpha_{i}}^{+} x_{\alpha_{i-1}}^{+}:$ or an operator $x_{\alpha_{i}}^{+}$
or an operator $x_{\alpha_{i-1}}^{+}$ or an identity
 is applied on every tensor component of $v_{\Lambda}$. For any $l,j=1,2,\ldots,n$ we have
 \begin{equation}\label{eqi:2}
 \left(E_{+}^{+}(-a_j,z)\otimes e^{\alpha_{j}}z^{\alpha_{j}}\right)v_{\Lambda_{l}}=
\begin{cases}
1\otimes e^{\alpha_{j}}v_{\Lambda_{l}}z, & \text{if }j=l;\\
1\otimes e^{\alpha_{j}}v_{\Lambda_{l}}, & \text{if }j\neq l.
\end{cases}
 \end{equation}
 Furthermore, all  the operators $E_{-}^{+}(-a_j,z)$ and $\phi_{j}(z)$, $j=1,2,\ldots,n$, consist of nonnegative powers of $z$.
 Considering the preceding observation, as well as  (\ref{eqi:1}) and (\ref{eqi:2}), we conclude that the lowest power of a variable $z_1$ 
 in (\ref{norm:exp}) equals $\sum_{s=1}^{m}\delta_{ij_{s}}$ and that the lowest power of a variable $z_2$ 
 in (\ref{norm:exp}) equals $\sum_{s=1}^{k}\delta_{i-1\hspace{1pt}j_{s}}$.
\end{dok}

Although we are mainly interested in finding relations among type $2$ quasi-particles, we state a simple corollary,
of a proof of Lemma \ref{tech5}, for type $1$ quasi-particles.

\begin{kor}
Let $\Lambda=c_{0}\Lambda_{0}+c_{j}\Lambda_{j}$ be a  highest weight as in  (\ref{weight}) and let $v_{\Lambda}$ be a highest weight vector  
of $L(\Lambda)$. We have
\begin{align*}
z_{1}^{\min\left\{m,k\right\}}B(z_2/z_1)x_{m\alpha_i}^{+}(z_{1})x_{k\alpha_{i-1}}^{+}(z_{2})v_{\Lambda}
\in z_{1}^{\sum_{s=1}^{m}\delta_{ij_{s}}}z_{2}^{\sum_{s=1}^{k}\delta_{i-1\hspace{1pt}j_{s}}}W(\Lambda)[[z_{1},z_{2}]].
\end{align*}
\end{kor}

The next lemma follows from Lemma \ref{tech5} and a fact that the operators $\bar{x}_{\alpha_{i}}^{+}(z_1)$ and  
$\bar{x}_{\alpha_{j}}^{+}(z_2)$, such that $a_{ij}=0,2$, commute (see Theorem \ref{106_ding_feigin}).
Its statement is an analogue of lemma 5.1 from \cite{G},  proven for affine Lie algebras of type $A_{n}^{(1)}$.

\begin{lem}\label{tech6}
Let $\Lambda=c_{0}\Lambda_{0}+c_{j}\Lambda_{j}$ be a highest weight as in  (\ref{weight}) and let $v_{\Lambda}$ be a highest weight vector  
of $L(\Lambda)$. For every  operator
$$\bar{x}^{+}_{m_{r_{n}^{(1)},n}\alpha_{n}}(z_{r_{n}^{(1)},n})\cdots \bar{x}^{+}_{m_{1,1}\alpha_{1}}(z_{1,1})$$
of color-charge-type
$$(m_{r_{n}^{(1)},n},\ldots,m_{1,n};\ldots;m_{r_{1}^{(1)},1},\ldots,m_{1,1})$$
and color-dual-charge-type
$$(r_{n}^{(1)},\ldots,r_{n}^{(k)};\ldots;r_{1}^{(1)},\ldots,r_{1}^{(k)})$$
there exist Taylor series 
$A_{i}^{(r,s)}(z)\in \mathbb{C}(q^{1/2})[[z]]$, where $i=2,3,..,n$, $r=1,2,\ldots,r_{i}^{(1)}$, $s=1,\ldots,r_{i-1}^{(1)}$,
such that $A_{i}^{(r,s)}(0)=1$ and
\begin{align}\label{the:exp}
&A(z_{r_{n}^{(1)},n},\ldots,z_{1,1})\left(            \prod_{i=2}^{n}   \prod_{r=1}^{r_{i}^{(1)}}        \prod_{s=1}^{r_{i-1}^{(1)}}       \prod_{t=1}^{\min\left\{m_{r,i},m_{s,i-1}\right\}}        \left(    1-q^{   2(  t+m_{r,i}-\min\left\{m_{r,i},m_{s,i-1}\right\}   )   }\frac{z_{s,i-1}}{z_{r,i}}    \right)              \right) \\ 
&\hspace{20pt}\cdot\bar{x}^{+}_{m_{r_{n}^{(1)},n}\alpha_{n}}(z_{r_{n}^{(1)},n})\cdots \bar{x}^{+}_{m_{1,1}\alpha_{1}}(z_{1,1})v_{\Lambda}\nonumber\\
&\hspace{40pt}\in\left(       \prod_{i=1}^{n}        \prod_{r=1}^{r_{i}^{(1)}}        z_{r,i}^{\sum_{s=1}^{m_{r,i}}\delta_{ij_{s}}   -\sum_{s=1}^{r_{i-1}^{(1)}}\min\left\{m_{r,i},m_{s,i-1}\right\}}              \right)W(\Lambda)[[z_{r_{n}^{(1)},n},\ldots,z_{1,1}]],\nonumber
\end{align}
where the Laurent series $A$ is given by
$$A(z_{r_{n}^{(1)},n},\ldots,z_{1,1}):=\prod_{i=2}^{n}   \prod_{r=1}^{r_{i}^{(1)}}        \prod_{s=1}^{r_{i-1}^{(1)}}A_{i}^{(r,s)}\left(z_{s,i-1}/z_{r,i}\right).$$
In the above formula (\ref{the:exp}) we   assume
$r_{0}^{(1)}:=0$ and $\sum_{s=1}^{0}\min\left\{m_{r,i},m_{s,i-1}\right\}:=0$.
\end{lem}


\subsection{Relations among quasi-particles of the same color}
Fix a color $i=1,2,\ldots,n$ and fix an integral dominant weight $\Lambda\in\hat{P}^{+}$.

\begin{lem}\label{relation_sc}
For any positive integers $m$ and $k$, $m\leq k$, the following $2m$ relations hold on every $L(\Lambda)$:
\begin{align*}
&\textrm{\normalfont(1)}&  \bar{x}_{m\alpha_i}^{+}(zq^{-2m})\bar{x}_{k\alpha_i}^{+}(z)=&\bar{x}_{(m+k)\alpha_i}^{+}(zq^{-2m}),\\
&\textrm{\normalfont(2)} & \bar{x}_{m\alpha_i}^{+}(zq^{-2(m-1)})\bar{x}_{k\alpha_i}^{+}(z)=&\bar{x}_{\alpha_i}^{+}(z)\bar{x}_{(m+k-1)\alpha_i}^{+}(zq^{-2(m-1)}),\\
&\hspace{6pt}\vdots&\vdots\\
&\textrm{\normalfont(m)}  &\bar{x}_{m\alpha_i}^{+}(zq^{-2(m-(m-1))})\bar{x}_{k\alpha_i}^{+}(z)=&\bar{x}_{(m-1)\alpha_i}^{+}(z)\bar{x}_{(k+1)\alpha_i}^{+}(zq^{-2(m-(m-1))}),\\
&\textrm{\normalfont(m+1)} & \bar{x}_{m\alpha_i}^{+}(zq^{2k})\bar{x}_{k\alpha_i}^{+}(z)=&\bar{x}_{(m+k)\alpha_i}^{+}(z),\\
&\textrm{\normalfont(m+2)}  &\bar{x}_{m\alpha_i}^{+}(zq^{2(k-1)})\bar{x}_{k\alpha_i}^{+}(z)=&\bar{x}_{\alpha_i}^{+}(zq^{2(k-1)})\bar{x}_{(m+k-1)\alpha_i}^{+}(z),\\
&\hspace{6pt}\vdots&\vdots\\
&\textrm{\normalfont(2m)} &\bar{x}_{m\alpha_i}^{+}(zq^{2(k-(m-1))})\bar{x}_{k\alpha_i}^{+}(z)=&\bar{x}_{(m-1)\alpha_i}^{+}(zq^{2(k-(m-1))})\bar{x}_{(k+1)\alpha_i}^{+}(z).
\end{align*}
\end{lem}

\begin{dok}
Relations (1)--(2m) of the lemma follow from the definition of the operator $\bar{x}^{+}_{m\alpha_{i}}(z)$ (see (\ref{type2})) and 
Corollary \ref{107_t2.comm.}. For example, we have
$$\bar{x}_{2\alpha_i}^{+}(zq^{-4})\bar{x}_{3\alpha_i}^{+}(z)=\bar{x}_{\alpha_i}^{+}(zq^{-4})\bar{x}_{\alpha_i}^{+}(zq^{-2})\bar{x}_{\alpha_i}^{+}(z)\bar{x}_{\alpha_i}^{+}(zq^{2})\bar{x}_{\alpha_i}^{+}(zq^{4})=\bar{x}_{5\alpha_i}^{+}(zq^{-4}),$$
which proves relation (1) when $m=2$ and $k=3$.
\end{dok}

 For every vector $v\in L(\Lambda)$ and for any  integers $m,k,N$ such that $1\leq m \leq k$
 we define a set
\begin{equation}
S_{N,v}^{m,k}:=\left\{\bar{x}_{m\alpha_{i}}^{+}(l)\bar{x}_{k\alpha_{i}}^{+}(N-l)v \hspace{4pt}|\hspace{4pt} l\in\mathbb{Z}\right\}.
\end{equation}

\begin{lem}\label{relations_sc2}
For any integers $m,k,N,r$ such that $1\leq m \leq k$ the vectors 
\begin{align*}
&\bar{x}_{m\alpha_{i}}^{+}(r)\bar{x}_{k\alpha_{i}}^{+}(N-r)v,\\
&\bar{x}_{m\alpha_{i}}^{+}(r+1)\bar{x}_{k\alpha_{i}}^{+}(N-(r+1))v,\\ 
&\hspace{80pt}\vdots\\
&\bar{x}_{m\alpha_{i}}^{+}(r+2m-1)\bar{x}_{k\alpha_{i}}^{+}(N-(r+2m-1))v
\end{align*} 
of the set $S_{N,v}^{m,k}$ can be expressed as a linear combination of some other vectors of the set $S_{N,v}^{m,k}$ and some vectors of the sets $S_{N,v}^{m',k'}$, 
such that 
$$0\leq m'<m,\quad m'+k'=m+k.$$ 
(For $m^{'}=0$ we set $\bar{x}_{0\alpha_{i}}^{+}(z):=1$.)
\end{lem}

\begin{dok}
Fix a vector $v$ and integers $m,k,N,r$.
First, we apply $2m$ relations of the Lemma \ref{relation_sc} on the vector $v$.
Next, we equate the coefficients of $z^{-N}$ in all  this relations, thus obtaining $2m$ equalities. Each of the
equalities consists of vectors 
\begin{equation}\label{seq:0}
\bar{x}_{m\alpha_{i}}^{+}(r)\bar{x}_{k\alpha_{i}}^{+}(N-r)v, \ldots, 
\bar{x}_{m\alpha_{i}}^{+}(r+2m-1)\bar{x}_{k\alpha_{i}}^{+}(N-(r+2m-1))v
\end{equation} 
and of other elements of the sets $S_{N,v}^{m,k}$ and $S_{N,v}^{m',k'}$,
where $0\leq m'<m$ and $m'+k'=m+k$. We can consider this equalities as a system of $2m$ linear equations in $2m$
variables (\ref{seq:0}). Since coefficients of variables are nonzero, we can assume that coefficients of the variable
$\bar{x}_{m\alpha_{i}}^{+}(r)\bar{x}_{k\alpha_{i}}^{+}(N-r)v$ in all  $2m$ equations are equal to $1$. Then the coefficient 
matrix of this system is a Vandermonde matrix whose determinant equals
$$(-1)^{m}\prod_{\substack{ r,s\in J\\ r< s}}\left(q^{2s}-q^{2r}\right),\textrm{ where }J:=\left\{1,2,\ldots,m, -k, -k+1,\ldots, -k+m-1\right\}.$$
Finally, we conclude that the matrix is regular, thus proving the lemma.
\end{dok}

The following corollary is an easy consequence of the above lemma.

\begin{kor}\label{relations_sc3}
For every vector $v\in L(\Lambda)$ and for any two integers $N$ and $m$, $m\geq 1$, 
 the vector 
 \begin{equation}\label{phd1}
 \bar{x}_{m\alpha_{i}}^{+}(l)\bar{x}_{m\alpha_{i}}^{+}(N-l)v,\qquad \textrm{where } N-l-2m< l \leq  N-l,    
 \end{equation}
 can be expressed as a linear combination of  vectors 
\begin{equation}\label{phd6}
\bar{x}_{m\alpha_{i}}^{+}(s)\bar{x}_{m\alpha_{i}}^{+}(N-s)v,\qquad \textrm{where }s\leq N-s-2m,
\end{equation}
 and some  vectors of the sets $S_{N,v}^{m',k'}$ 
such that
$$0\leq m'<m,\quad m'+k'=2m.$$ 
\end{kor}


\subsection{The spanning set of \texorpdfstring{$W(\Lambda)$}{W(Lambda)}}
The following lemma gives us two important properies of the orders defined at the end of Subsection \ref{sub:order}.

\begin{lem}\label{boundary}
For every monomial
$b\in\bar{\mathfrak{S}}_{W(\Lambda)}^{(c)}$
there are
\begin{enumerate}
  \item finitely many monomials $b^{'}\in\bar{\mathfrak{S}}_{W(\Lambda)}^{(c)}$ of the same color-charge-type and the same degree 
  as a monomial  $b$, such that
$$b\prec b^{'}\qquad\textrm{and}\qquad b^{'}v_{\Lambda}\neq 0\textrm{;}$$
\item finitely many color-charge-types of the same color-type as a monomial $b$, that are greater, regarding partial order ``$\prec$'' or linear order ``$<$'', than color-charge-type of a monomial $b$.
\end{enumerate}
\end{lem}

\begin{dok}
(1) Let $b=b_{1}b_{2}\cdots b_{n}\in \bar{\mathfrak{S}}_{W(\Lambda)}^{(c)}$ ($b^{'}=b_{n}^{'}\cdots b_{1}^{'}\in \bar{\mathfrak{S}}_{W(\Lambda)}^{(c)}$) 
be a monomial written as a product of monochrome monomials $b_{i}$ ($b_{i}^{'}$)
of color $i=1,2,\ldots,n$. Denote by $l_{i}$ ($l_{i}^{'}$) a degree of the monomial $b_i$ ($b_{i}^{'}$). 
Suppose $b$ and $b^{'}$ are of the same color-charge-type and $b\prec b^{'}$.
Obviously, $l_1\leq l_{1}^{'}$. The integer $l_{1}^{'}$ is also bounded from above (see Lemma \ref{restricted}) so it can
have only  finitely many integer values. For every  $l_{1}^{'}$ there exist finitely many $l_{1}^{'}$ degree monomials $b_{1}^{'}$ of the same 
color-charge-type as $b_1$, such that $b_{1}^{'}v_{\Lambda}\neq 0$.

We can analogously prove  that for every monomial $b_{1}^{'}$ of color $1$ there exist finitely many monomials $b_{2}^{'}$ 
of the same color-charge-type as $b_{2}$ and such that $b_{2}b_{1}\prec b_{2}^{'}b_{1}^{'}$
and, more generally, that for every monomial
$b_{i}^{'}\cdots b_{1}^{'}$
of the same color-charge-type as $b_{i}\cdots b_{1}$, such that
 $b_{i}\cdots b_{1}\prec b_{i}^{'}\cdots b_{1}^{'}$,
there exist finitely many monomials $b_{i+1}^{'}$ 
of the same color-charge-type as $b_{i+1}$ such that $b_{i+1}\cdots b_{1}\prec b_{i+1}^{'}\cdots b_{1}^{'}$. The proof follows by induction.

(2) The second statement follows from the fact that every positive integer 
has (only) finite number of partitions. 
\end{dok}

\begin{thm}\label{spanning_set}
For a given highest weight $\Lambda=c_{0}\Lambda_{0}+c_{j}\Lambda_{j}$ as in (\ref{weight}) the set 
$$\left\{bv_{\Lambda}\hspace{2pt}\left|\hspace{2pt}\right. b\in\bar{\mathfrak{B}}_{W(\Lambda)}\right\}$$ 
spans  the principal subspace $W(\Lambda)$. 
\end{thm}

\begin{dok}
Lemmas \ref{tech6} and \ref{relations_sc2} as well as Corollary \ref{relations_sc3} will allow us
to prove the theorem in the same way  Georgiev proved Theorem 5.1 in \cite{G}. 
Since the principal subspace $W(\Lambda)$ is a direct sum of its weight subspaces,
$$W(\Lambda)=\bigoplus_{\mu\in \hat{P}}W(\Lambda)_\mu,
\quad W(\Lambda)_\mu=\left\{v\in W(\Lambda)\hspace{2pt}\left|\hspace{2pt}\right.q^{h}v=q^{\mu(h)}v\textrm{ za sve }h\in\hat{P}^{\vee}\right\},$$ 
it is sufficient to prove that every homogenous vector $v\in W(\Lambda)$ is a linear combination of some vectors 
$bv_{\Lambda}$, $b\in\bar{\mathfrak{B}}_{W(\Lambda)}$. The proof will be carried out by an induction on color-charge-types of monomials 
$b\in\bar{\mathfrak{S}}_{W(\Lambda)}^{(c)}$ (ordered by ``$\prec$'').

(I) Let $b$ be a monomial in $\bar{\mathfrak{S}}_{W(\Lambda)}^{(c)}$ containing a quasi-particle of color $i$, charge $m_{r,i}$
and of degree $l_{r,i}$ that violates a condition
\begin{equation}\label{some:cond}
l_{r,i}\leq \sum_{s=1}^{r_{i-1}^{(1)}}\min\left\{m_{r,i},m_{s,i-1}\right\}-\sum_{s=1}^{m_{r,i}}\delta_{ij_{s}}-m_{r,i}.
\end{equation}
Suppose $b$ is of a color-charge-type 
$(m_{r_{n}^{(1)},n},\ldots,m_{1,n};\ldots;m_{r_{1}^{(1)},1},\ldots,m_{1,1})$
and suppose $(l_{r_{n}^{(1)},n},\ldots,l_{1,n};\ldots;l_{r_{1}^{(1)},1},\ldots,l_{1,1})$ is sequence of 
degrees of its quasi-particles. Recall Lemma \ref{tech6}.
A coefficient of
\begin{equation}\label{some:var}
z_{r_{n}^{(1)},n}^{-l_{r_{n}^{(1)},n}-m_{r_{n}^{(1)},n}}\cdots z_{1,n}^{-l_{1,n}-m_{1,n}}\cdots z_{r_{1}^{(1)},1}^{-l_{r_{1}^{(1)},1}-m_{r_{1}^{(1)},1}}\cdots z_{1,1}^{-l_{1,1}-m_{1,1}}.
\end{equation}
in (\ref{the:exp}) equals a certain linear combination of monomials acting on the vector $v_{\Lambda}$. 
One of them is  the vector $bv_{\Lambda}$ while the others, $b^{'}v_{\Lambda}$, satisfy $b\prec b^{'}$. 
Note that the relation $b\prec b^{'}$ is actually an easy consequence of the following fact:
quotients of variables $z/w$ appearing in (\ref{the:exp}) are such that the operator in variable $z$
stands on the right side of the operator in variable $w$.
Since the monomial $b$ does not satisfy the condition (\ref{some:cond}), Lemma \ref{tech6} implies
the above-mentioned linear combination (i.e. the coefficient of (\ref{some:var})) being equal to $0$. 
Therefore, we can express the vector $bv_{\Lambda}$ as a linear combination of vectors 
$b^{'}v_{\Lambda}$ such that $b\prec b^{'}$ and such that the monomials 
$b^{'}\in \bar{\mathfrak{S}}_{W(\Lambda)}^{(c)}$ are of the same color-charge-type and the same degree  
as the monomial $b$ (but different color-degree-type).

(II) Suppose $b\in \bar{\mathfrak{S}}_{W(\Lambda)}^{(c)}$ satisfies (\ref{some:cond}) but contains a quasi-particle
$\bar{x}_{m_{r,i}\alpha_{i}}^{+}(l_{r,i})$ that violates a condition
\begin{equation}\label{some:cond2}
l_{r,i}\leq \sum_{s=1}^{r_{i-1}^{(1)}}\min\left\{m_{r,i},m_{s,i-1}\right\}-\sum_{s=1}^{m_{r,i}}\delta_{ij_{s}}-\sum_{m_{t,i}>m_{r,i}}2m_{r,i}-m_{r,i}.
\end{equation}
By applying Lemma \ref{relations_sc2} on all pairs consisting of the quasi-particle $\bar{x}_{m_{r,i}\alpha_{i}}^{+}(l_{r,i})$
and some other quasi-particle of color $i$ and charge greater than $m_{r,i}$ we can write a vector $bv_{\Lambda}$
as a linear combination of
\begin{enumerate}
  \item vectors $b^{'}v_{\Lambda}$ whose monomials $b^{'}$ satisfy (\ref{some:cond2}) and have 
   the same color-charge-type and the same color-degree-type as $b$; 
  \item vectors $b^{''}v_{\Lambda}$ whose monomials $b^{''}$ satisfy $b\prec b^{''}$ and have the same color-degree-type 
  and the same color-type as $b$
  but they are not of the same color-charge-type as $b$;
  \item vectors $b^{'''}v_{\Lambda}$ whose monomials  $b^{'''}$ do not satisfy (\ref{some:cond2}) 
 but have the same color-charge-type and the same color-degree-type as $b$.
\end{enumerate}
Now we only have to explain how to deal with the vectors $b^{'''}v_{\Lambda}$. 
Notice that we could have applied Lemma \ref{relations_sc2} in such a way  that the obtained monomials $b^{'''}$
not only violate condition (\ref{some:cond2}) but also  condition (\ref{some:cond}). Therefore, we can 
apply step (I) on the vectors $b^{'''}v_{\Lambda}$, thus writing them as a linear combination of some new 
vectors $b^{''''}v_{\Lambda}$
whose monomials $b^{''''}$  satisfy $b^{'''}\prec b^{''''}$ and  $b\prec b^{''''}$.

(III) Suppose that the monomial $b=\ldots\bar{x}^{+}_{m_{r+1,i}\alpha_{i}}(l_{r+1,i})\bar{x}^{+}_{m_{r,i}\alpha_{i}}(l_{r,i})\ldots$
does not satisfy
\begin{equation}\label{some:cond3}
l_{r+1,i}\leq l_{r,i}-2m_{r,i}.
\end{equation}
By applying Corollary \ref{relations_sc3} we can write a vector $bv_{\Lambda}$
as a linear combination of
\begin{enumerate}
  \item vectors $b^{'}v_{\Lambda}$ whose monomials $b^{'}$ satisfy (\ref{some:cond3}) and have 
   the same color-charge-type and the same color-degree-type as $b$; 
   \item vectors $b^{''}v_{\Lambda}$ whose monomials $b^{'}$ satisfy $b\prec b^{''}$ and have
   the same color-degree-type as $b$.
\end{enumerate}

Notice that the properties of the order ``$\prec$'' given by Lemma 
\ref{boundary} guarantee that, after finitely many steps (I)--(III), we can write the vector
$bv_{\Lambda}$, $b\in \bar{\mathfrak{S}}_{W(\Lambda)}^{(c)}$, as a linear combination of  vectors
of a set $\left\{bv_{\Lambda}\hspace{2pt}\left|\hspace{2pt}\right. b\in\bar{\mathfrak{B}}_{W(\Lambda)}\right\}$.
Since the vectors $bv_{\Lambda}$, $b\in \bar{\mathfrak{S}}_{W(\Lambda)}^{(c)}$, span the principal subspace $W(\Lambda)$ 
(Lemma \ref{simple_lemma2}), the statement of the theorem follows.
\end{dok}


\section{Linear independence of the set \texorpdfstring{$\mathfrak{B}_{W(\Lambda)}$}{BW(Lambda)}}


\subsection{Projection \texorpdfstring{$\pi$}{pi}}
For every dominant integral highest weight $\Lambda=c_{0}\Lambda_{0}+c_{j}\Lambda_{j}$ as in (\ref{weight}) 
the principal subspace $W(\Lambda)$ has a realization as a subspace of the tensor product of $c$ level one principal subspaces
\begin{equation*}
W(\Lambda)\subset W(\Lambda_{j_{1}})\otimes \ldots\otimes W(\Lambda_{j_{c}}),
\end{equation*}
where
$$v_\Lambda= v_{\Lambda_{j_{1}}}\otimes \ldots\otimes v_{\Lambda_{j_{c}}}$$
and indices $j_{s}$ are defined by (\ref{indices}).
Consider the direct sum decomposition
 $$W(\Lambda_{j_{1}})\otimes \ldots\otimes W(\Lambda_{j_{c}})=
 \bigoplus_{  \substack{   r_{n}^{(1)},\ldots, r_{1}^{(1)}\geq 0\\ \ldots\\r_{n}^{(c)},\ldots, r_{1}^{(c)}\geq 0    }   }
  W(\Lambda_{j_{1}})_{(r_{n}^{(1)},\ldots, r_{1}^{(1)})}\otimes \cdots\otimes W(\Lambda_{j_{c}})_{(r_{n}^{(c)},\ldots, r_{1}^{(c)})},$$
  where
    $W(\Lambda_{j_{s}})_{(r_{n}^{(s)},\ldots, r_{1}^{(s)})}:=W(\Lambda_{j_{s}})_{\Lambda_{j_{s}}+\sum_{i=1}^{n}r_{i}^{(s)}\alpha_{i}}$
  is the weight subspace
   $$W(\Lambda_{j_{s}})_{\Lambda_{j_{s}}+\sum_{i=1}^{n}r_{i}^{(s)}\alpha_{i}}:=\left\{v\in W(\Lambda_{j_{s}})\hspace{2pt}|\hspace{2pt} K_{i}v=q^{(\Lambda_{j_{s}}+\sum_{i=1}^{n}r_{i}^{(s)}\alpha_{i})(\alpha_{i}^{\vee})}v\textrm{ for }i=1,2,\ldots,n\right\}.$$ 
For every color-dual-charge-type 
$(r_{n}^{(1)},\ldots,r_{n}^{(c)};\ldots;r_{1}^{(1)},\ldots,r_{1}^{(c)})$
the decomposition above gives us a projection
\begin{equation*}
 \pi_{(r_{n}^{(1)},\ldots,r_{1}^{(c)})}\colon W(\Lambda_{j_{1}})\otimes \ldots\otimes W(\Lambda_{j_{c}})\to
  W(\Lambda_{j_{1}})_{(r_{n}^{(1)},\ldots, r_{1}^{(1)})}\otimes \cdots\otimes W(\Lambda_{j_{c}})_{(r_{n}^{(c)},\ldots, r_{1}^{(c)})}.
 \end{equation*}
 The projection can be in an obvious way generalized to the space of formal Laurent series with coefficients in $W(\Lambda_{j_{1}})\otimes \ldots\otimes W(\Lambda_{j_{c}})$.
  For an operator 
\begin{equation}\label{5:op}
x^{+}_{m_{r_{n}^{(1)},n}\alpha_{n}}(z_{r_{n}^{(1)},n})\cdots x^{+}_{m_{1,1}\alpha_{1}}(z_{1,1})
\end{equation}
of color-dual-charge-type
$(r_{n}^{(1)},\ldots,r_{n}^{(c)};\ldots;r_{1}^{(1)},\ldots,r_{1}^{(c)})$
we can, by using Lemma \ref{x_order}, write down the action of the projection $\pi_{(r_{n}^{(1)},\ldots,r_{1}^{(c)})}$
on the series
 $$x^{+}_{m_{r_{n}^{(1)},n}\alpha_{n}}(z_{r_{n}^{(1)},n})\cdots x^{+}_{m_{1,1}\alpha_{1}}(z_{1,1})v_{\Lambda}=x^{+}_{m_{r_{n}^{(1)},n}\alpha_{n}}(z_{r_{n}^{(1)},n})\cdots x^{+}_{m_{1,1}\alpha_{1}}(z_{1,1})v_{\Lambda_{j_1}}\otimes\ldots\otimes v_{\Lambda_{j_c}}.$$
First, choose an operator in (\ref{5:op}), for example $x^{+}_{m_{l,i}\alpha_{i}}(z_{l,i})$, where $i=1,2,\ldots,n$, $l=1,2,\ldots, r_{i}^{(1)}$.
Recall Definition \ref{type_1_op.}.
In the formula
$$\pi_{(r_{n}^{(1)},\ldots,r_{1}^{(c)})}\left(x^{+}_{m_{r_{n}^{(1)},n}\alpha_{n}}(z_{r_{n}^{(1)},n})\cdots x^{+}_{m_{1,1}\alpha_{1}}(z_{1,1})v_{\Lambda_{j_1}}\otimes\ldots\otimes v_{\Lambda_{j_c}}\right)$$
the operator $x^{+}_{m\alpha_{i}}(z):=x^{+}_{m_{l,i}\alpha_{i}}(z_{l,i})$ will correspond to the term
\begin{align}\label{5:for}
 &\hspace{-20pt}\lim_{z_{p}\to zq^{2(p-1)}}\Bigg(\Bigg.\left(\prod_{r=1}^{m-1}\prod_{s=r+1}^{m}
\left(1-q^{2}\frac{z_{s}}{z_{r}}\right)\right)\\
&\hspace{11pt}\cdot\phi_{i}(z_{1}q^{\frac{1}{2}})\phi_{i}(z_{2}q^{\frac{1}{2}})\cdots\phi_{i}(z_{m-1}q^{\frac{1}{2}})x_{\alpha_{i}}^{+}(z_{m})\nonumber\\
&\hspace{5pt}\otimes\phi_{i}(z_{1}q^{\frac{3}{2}})\phi_{i}(z_{2}q^{\frac{3}{2}})\cdots\phi_{i}(z_{m-2}q^{\frac{3}{2}})x_{\alpha_{i}}^{+}(z_{m-1}q)\nonumber\\
&\hspace{130pt}\vdots\nonumber\\
&\hspace{5pt}\otimes\phi_{i}(z_{1}q^{m-\frac{3}{2}})x_{\alpha_{i}}^{+}(z_{2}q^{m-2})\nonumber\\
&\hspace{5pt}\otimes x_{\alpha_{i}}^{+}(z_{1}q^{m-1})\nonumber\\
&\hspace{5pt}\otimes \underbrace{1\otimes\cdots\otimes 1}_{\mbox{$c-m$}} \Bigg.\Bigg)\nonumber.
 \end{align}
Notice that the projection $\pi_{(r_{n}^{(1)},\ldots,r_{1}^{(c)})}$ forces the operators 
 $x_{\alpha_{i}}^{+}(z_j)$, $j=1,2,\ldots,m=m_{l,i}$, to spread along the 
 $m$ leftmost tensor factors in an order established by Lemma \ref{x_order}. 
 Each of these $m$ tensor factors contains exactly one operator $x_{\alpha_{i}}^{+}(z_j)$.
 
 The projection was defined in a similar way as in \cite{G}. The only difference is a fact that by applying the original projection 
 the operators $x_{\alpha_{i}}^{+}(z_j)$ in formula \ref{5:for}  spread along the $m$ rightmost tensor factors. 
 This small modification will allow us to carry out the linear independence proof at the end of this section, although
 the Hopf algebra structure on  $U_{q}(\widehat{\mathfrak{sl}}_{n+1})$ is somewhat more complicated then the one on
 $U(\widehat{\mathfrak{sl}}_{n+1})$.


\subsection{Operator \texorpdfstring{$\mathcal{Y}$}{Y}}
In \cite{Koyama} Y. Koyama found a realization of vertex operators for  level one integrable
highest weight modules of $U_{q}(\widehat{\mathfrak{sl}}_{n+1})$. We will consider here a similar operator defined
on a space 
$$W:=K(1)\otimes \mathbb{C}\left\{P\right\}.$$
The main properties of this operator will be relations given by Theorem \ref{relations_Y}.  
They will allow us to use the operator in the linear independence proof in the next subsection.
 The proofs of Lemma \ref{5:lemx} and Theorem \ref{relations_Y} follow from a direct calculation and they are, therefore, omitted.

Let $i=1,2,\ldots,n$ and $l\in\mathbb{Z}$, $l\neq 0$. Define elements $a_{i}^{*}(l)\in U_{q}(\hat{\mathfrak{h}})$ by
\begin{align*}
a_{i}^{*}(l):=&\frac{[l][(n-i+1)l]}{[(n+1)l][l]}a_{1}(l)+\frac{[2l][(n-i+1)l]}{[(n+1)l][l]}a_{2}(l)+\ldots\\
&\hspace{20pt}\ldots+\frac{[(i-1)l][(n-i+1)l]}{[(n+1)l][l]}a_{i-1}(l)+\frac{[il][(n-i+1)l]}{[(n+1)l][l]}a_{i}(l)\\
&\hspace{40pt}+\frac{[il][(n-i)l]}{[(n+1)l][l]}a_{i+1}(l)+\ldots+\frac{[il][l]}{[(n+1)l][l]}a_{n}(l).
\end{align*}

\begin{lem}\label{5:lemx}
For any colors $i,j=1,2,\ldots,n$ and integers $l,k$ 
\begin{equation*}
[a_{i}^{*}(l),a_{j}(k)]=\delta_{ij}\delta_{l+k\hspace{2pt}0}\frac{[l]^{2}}{l}.
\end{equation*}
\end{lem}

Fix a color $i=1,2,\ldots,n$.  We define the following operators on the space $W$:
\begin{align*}
E_{-}(a_{i}^{*},z)&:=\exp\left(\sum_{r=1}^{\infty}\frac{q^{r/2}}{[r]}a_{i}^{*}(-r)z^{r}\right),\\
E_{+}(a_{i}^{*},z)&:=\exp\left(-\sum_{r=1}^{\infty}\frac{q^{r/2}}{[r]}a_{i}^{*}(r)z^{-r}\right).
\end{align*}
Denote by $W\left\{z\right\}$ the space
$$W\left\{z\right\}:=\left\{\sum_{h\in\mathbb{C}}v_{h}z^{h}\hspace{4pt}\Big|\hspace{4pt} v_{h}\in W\textrm{ for all }h\in\mathbb{C}\right\}.$$

\begin{defn}
We define an operator 
$\mathcal{Y}(e^{\lambda_{i}},z)\in \om(W,W\left\{Z\right\})$
by
\begin{align*}
\mathcal{Y}(e^{\lambda_{i}},z):=E_{-}(a_{i}^{*},z)E_{+}(a_{i}^{*},z)
\otimes e^{\lambda_{i}}(-1)^{(1-\delta_{in})i \partial_{\lambda_{n}}}z^{\partial_{\lambda_{i}}}.
\end{align*}
\end{defn}

Notice that for every vector $w\in W$ the series $\mathcal{Y}(e^{\lambda_{i}},z)w$ contains a finite number of negative powers of variable $z$.
The following theorem can be proved by a direct calculation.

\begin{thm}\label{relations_Y}
For any colors $i,j=1,2,\ldots,n$ the following relations hold on $W$:
\begin{enumerate}
  \item $\displaystyle[x_{\alpha_{i}}^{+}(z_{1}),\mathcal{Y}(e^{\lambda_{j}},z_2)]=0,$
  \item $[x_{\alpha_{i}}^{-}(z_{1}),\mathcal{Y}(e^{\lambda_{j}},z_2)]=0\textrm{ if }i\neq j,$
  \item $(z_{1}-qz_{2})x_{\alpha_{i}}^{-}(z_{1})\mathcal{Y}(e^{\lambda_{i}},z_{2})=(qz_{1}-z_{2})\mathcal{Y}(e^{\lambda_{i}},z_{2})x_{\alpha_{i}}^{-}(z_{1}),$
  \item $[\phi_{i}(z_{1}),\mathcal{Y}(e^{\lambda_{j}},z_2)]=[\psi_{i}(z_{1}),\mathcal{Y}(e^{\lambda_{j}},z_2)]=0\textrm{ if }i\neq j,$
  \item $(q^{1/2} z_{1}-qz_{2})\phi_{i}(z_1)\mathcal{Y}(e^{\lambda_{i}},z_2)=(q^{3/2}z_{1}-z_{2})\mathcal{Y}(e^{\lambda_{i}},z_2)\phi_{i}(z_1),$
  \item $(z_{1}-q^{3/2}z_{2})\psi_{i}(z_1)\mathcal{Y}(e^{\lambda_{i}},z_2)=(qz_{1}-q^{1/2}z_{2})\mathcal{Y}(e^{\lambda_{i}},z_2)\psi_{i}(z_1).$
\end{enumerate}
\end{thm}


\subsection{Proof of linear independence} 
We begin with a list of relations we will use in the proof of a linear independence (theorem \ref{independent}).

\begin{lem}
On every level one integrable highest weight module we have:
\begin{align}
&x^{+}_{\alpha_i}(r)e^{\lambda_j}=\varepsilon_{ij}e^{\lambda_j}x^{+}_{\alpha_i}(r+\delta_{ij});\label{inv:1}\\
&x^{+}_{\alpha_i}(r)e^{\alpha_j}=(-1)^{(\alpha_{i},\alpha_{j})}e^{\alpha_j}x^{+}_{\alpha_i}(r+2\delta_{ij}-\delta_{i\hspace{1pt}j-1}-\delta_{i\hspace{1pt}j+1});\label{inv:2}\\
&\phi_{i}(s)e^{\lambda_j}=q^{-\delta_{ij}}e^{\lambda_j}\phi_{i}(s);\label{inv:3}\\
&\phi_{i}(s)e^{\alpha_j}=q^{-2\delta_{ij}+\delta_{i\hspace{1pt}j-1}+\delta_{i\hspace{1pt}j+1}}e^{\alpha_j}\phi_{i}(s)\label{inv:4}
\end{align}
for some $\varepsilon_{ij}=\pm 1$ and for all $i,j=1,2,\ldots,n$, $r\in\mathbb{Z}$, $s\in\mathbb{Z}_{\leq 0}$.
\end{lem}
All of the above relations can be proved by a simple calculation.

\begin{thm}\label{independent}
For a given highest weight $\Lambda=c_{0}\Lambda_{0}+c_{j}\Lambda_{j}$ as in (\ref{weight}) the set 
$$\left\{bv_{\Lambda}\hspace{2pt}\left|\hspace{2pt}\right. b\in\mathfrak{B}_{W(\Lambda)}\right\}$$ 
is linearly independant. 
\end{thm}

\begin{dok}
Relations (1) and (4) of  Theorem \ref{relations_Y}, as well as the projection $\pi_{(r_{n}^{(1)},\ldots,r_{1}^{(c)})}$, allow us to 
carry out the proof in the same way  Georgiev proved Theorem 5.2 in \cite{G}. 

Let $b\in \mathfrak{B}_{W(\Lambda)}$ be a monomial 
$$b=x^{+}_{m_{r_{n}^{(1)},n}\alpha_{n}}(l_{r_{n}^{(1)},n})\cdots x^{+}_{m_{1,n}\alpha_{n}}(l_{1,n})\cdots x^{+}_{m_{r_{1}^{(1)},1}\alpha_{1}}(l_{r_{1}^{(1)},1})\cdots x^{+}_{m_{1,1}\alpha_{1}}(l_{1,1})$$
of color-charge-type
$$(m_{r_{n}^{(1)},n},\ldots,m_{1,n};\ldots;m_{r_{1}^{(1)},1},\ldots,m_{1,1})$$
and corresponding color-dual-charge type
$$(r_{n}^{(1)},\ldots,r_{n}^{(c)};\ldots;r_{1}^{(1)},\ldots,r_{1}^{(c)}).$$
First, we prove $bv_{\Lambda}\neq 0$.

Suppose $bv_{\Lambda}=0$. Then $\pi_{(r_{n}^{(1)},\ldots,r_{1}^{(c)})}bv_{\Lambda}=0$.
A positive integer  $m:=m_{1,1}$ is a maximal charge of the color $1$ quasi-particles  in the monomial $b$. 
Consider an action of 
\begin{equation}\label{Y}
\underbrace{1\otimes \ldots\otimes 1 }_\textrm{$m -1$}\otimes\rez_{z}\left(   z^{-1-(\lambda_{1},\lambda_{j_{m}}) }       \mathcal{Y}(e^{\lambda_{1}},z)        \right)\otimes \underbrace{1\otimes\cdots\otimes 1}_\textrm{$c-m-1$}
\end{equation} 
on $\pi_{(r_{n}^{(1)},\ldots,r_{1}^{(c)})}bv_{\Lambda}$.
The $m$-th tensor component of (\ref{Y}) commutes with all  the operators of the $m$-th component of
$\pi_{(r_{n}^{(1)},\ldots,r_{1}^{(c)})}bv_{\Lambda}=\pi_{(r_{n}^{(1)},\ldots,r_{1}^{(c)})}b(v_{\Lambda_{j_1}}\otimes \ldots\otimes v_{\Lambda_{j_c}})$
acting on the vector $v_{\Lambda_{j_m}}$ (see (1) in Theorem \ref{relations_Y}).
Therefore, we can move an operator 
$\rez_{z}\left(   z^{-1-(\lambda_{1},\lambda_{j_{m}}) }       \mathcal{Y}(e^{\lambda_{1}},z)        \right)$
all the way to the right. Notice that
$$\mathcal{Y}(e^{\lambda_{1}},z)v_{\Lambda_{j_m}}=Ce^{\lambda_{1}}v_{\Lambda_{j_m}},$$
where $C\in\mathbb{C}(q^{1/2})$ is a nonzero constant. By employing (\ref{inv:1}) and (\ref{inv:3}) we can move an  operator
$e^{\lambda_{1}}$ all the way to the left, thus getting
$$Ce^{\lambda_{1}}\pi_{(r_{n}^{(1)},\ldots,r_{1}^{(c)})}b^{'}v_{\Lambda}=0,$$
where $b^{'}\in \mathfrak{B}_{W(\Lambda)}$ is  obtained from the monomial $b$ by adding $1$ to the degrees of its color $1$ quasi-particles,
$$b^{'}=x^{+}_{m_{r_{n}^{(1)},n}\alpha_{n}}(l_{r_{n}^{(1)},n})\cdots  x^{+}_{m_{1,2}\alpha_{2}}(l_{1,2}) x^{+}_{m_{r_{1}^{(1)},1}\alpha_{1}}(l_{r_{1}^{(1)},1}+1)\cdots x^{+}_{m_{1,1}\alpha_{1}}(l_{1,1}+1).$$
Removing the invertible operator $e^{\lambda_{1}}$ and the constant $C\neq 0$ we get 
$$\pi_{(r_{n}^{(1)},\ldots,r_{1}^{(c)})}b^{'}v_{\Lambda}=0.$$ 

By repeating the above described algorithm we can, step by step, increase the degrees of all  the color $1$ quasi-particles in $b$.
Of course, in every step we get a new monomial that is an element of $\mathfrak{B}_{W(\Lambda)}$. 
We stop the algorithm  when the degree of the rightmost quasi-particle becomes equal to 
$-m-\sum_{r=1}^{m}\delta_{1j_{r}}$. Denote a corresponding monomial by $b^{''}$. Since
$x^{+}_{\alpha_{i}}(-1-\delta_{ij})v_{\Lambda_{j}}\neq 0$ and $x^{+}_{\alpha_{i}}(-\delta_{ij})v_{\Lambda_{j}}= 0$
for $i,j=1,2,\ldots,n$, an integer $-m-\sum_{r=1}^{m}\delta_{1j_{r}}$  is the maximal degree for which
 the corresponding quasi-particle does not annihilate
$v_{\Lambda}$. Dropping the rightmost quasi-particle of the monomial $b^{''}$ we get a monomial $b^{'''}$ of
color-charge-type
$$(m_{r_{n}^{(1)},n},\ldots,m_{1,n};\ldots;m_{r_{2}^{(1)},2},\ldots,m_{1,2};m_{r_{1}^{(1)},1},\ldots,m_{2,1})$$
and dual-color-charge-type
$$(r_{n}^{(1)},\ldots,r_{n}^{(c)};\ldots;r_{2}^{(1)},\ldots,r_{2}^{(c)};r_{1}^{(1)}-1,\ldots,r_{1}^{(c)}-1).$$
We have
\begin{align*}
0&=\pi_{(r_{n}^{(1)},\ldots,r_{1}^{(c)})}  b^{''}x^{+}_{m\alpha_{1}}(-m-\sum_{r=1}^{m}\delta_{1j_{r}})v_{\Lambda}\\
&=D\pi_{(r_{n}^{(1)},\ldots,r_{1}^{(c)})}  b^{''}\left(\underbrace{e^{\alpha_{1}}\otimes\cdots\otimes e^{\alpha_{1}}}_\textrm{$m$}\otimes 1\otimes\cdots\otimes 1\right)v_{\Lambda}
\end{align*}

for some nonzero constant $D$.
Formulas (\ref{inv:2}) and (\ref{inv:4}) allow us to move the operators $e^{\alpha_{1}}$ all the way to the left.
 Of course, by doing this we will change the degrees of the monomial $b^{'''}$ in the following way:
 \begin{itemize}
   \item the degree of any color $1$ quasi-particle will increase by a double value of its charge;
   \item the degree of any color $2$ quasi-particle wil decrease by a value of its charge;
   \item the degree of any color $i=3,4,\ldots,n$ quasi-particle will remain the same.
 \end{itemize}
 Denote a (new) monomial, that has the modified degrees, by $b^{''''}$ and then notice that $b^{''''}\in \mathfrak{B}_{W(\Lambda)}$.
 Dropping the invertible operators  $e^{\alpha_1}$ and a nonzero constant we get
 $$\pi_{(r_{n}^{(1)},\ldots ,r_{2}^{(c)},r_{1}^{(1)}-1,\ldots,r_{1}^{(c)}-1)}b^{''''}v_{\Lambda}=0.$$

By comparing the monomials $b$ and $b^{''''}$ we see that $b^{''''}$ lacks the rightmost quasi-particle of
$b$ and that $b^{''''}$ has somewhat modified degrees of its other quasi-particles.
In the same way we remove, step by step, all  the color $1$ quasi-particles of the monomial $b$. 
Then we remove all  the color $2$ quasi-particles of $b$ and so on. At the end, by removing all the color $n$
quasi-particles we get $v_{\Lambda}=0$. Contradiction! We conclude $bv_{\Lambda}\neq 0$ for $b\in \mathfrak{B}_{W(\Lambda)}$.

Now assume that 
\begin{equation}\label{lin.comb.}
\sum_{s=1}^{r}a_{s}b_{s}v_{\Lambda}=0
\end{equation}
for some $b_{1},\ldots,b_{r}\in\mathfrak{B}_{W(\Lambda)}$ and $a_{1},\ldots,a_{r}\in\mathbb{C}(q^{1/2})\setminus\left\{0\right\}$.
 Since $W(\Lambda)$ is a direct sum of its weight spaces, we can assume that the monomials $b_{s}$, $s=1,2,\ldots,r$,
 have the same color-type and the same degree. Suppose $b_{1}<b_{s}$ for $s=2,3,\ldots,r$. Now we can carry out the above described
 algorithm of charge reduction on (\ref{lin.comb.}) not stopping until we remove all  the charges of $b_{1}$. 
 By doing this the first summand $a_{1}b_{1}v_{\Lambda}$ is replaced by $Ca_{1}v_{\Lambda}$ for some $C\in\mathbb{C}(q^{1/2})$, $C\neq 0$.
 Notice that all the other monomials get annihilated at some intermediate stage of the charge reduction so (\ref{lin.comb.}) is replaced by
 $Ca_{1}v_{\Lambda}=0$.
 This implies $a_{1}=0$. Contradiction! The theorem now follows.
 \end{dok}

Now we can prove our main result.

\begin{thm}\label{main}
For a given highest weight $\Lambda=c_{0}\Lambda_{0}+c_{j}\Lambda_{j}$ as in (\ref{weight}) the sets 
$$\left\{bv_{\Lambda}\hspace{2pt}\left|\hspace{2pt}\right. b\in\mathfrak{B}_{W(\Lambda)}\right\}\qquad\textrm{and} \qquad \left\{bv_{\Lambda}\hspace{2pt}\left|\hspace{2pt}\right. b\in\bar{\mathfrak{B}}_{W(\Lambda)}\right\}$$ 
form the bases for the principal subspace $W(\Lambda)$. 
\end{thm}

\begin{dok}
The theorem is a consequence of Theorem \ref{spanning_set} and Theorem \ref{independent}. 
For a monomial $b\in\mathfrak{B}_{W(\Lambda)}$
denote by $\bar{b}$ a monomial in $\bar{\mathfrak{B}}_{W(\Lambda)}$ that 
has the same color-charge-type and the same degrees of its quasi-particles as $b$. 
Notice that the  vectors $bv_{\Lambda}$ and $\bar{b}v_{\Lambda}$
are weight vectors of the same weight. 
Since a principal subspace $W(\Lambda)\subset L(\Lambda)$ is a direct sum of finite dimensional weight subspaces, 
for every weight $\mu$ of $W(\Lambda)$ there are finitely many monomials $\bar{b}_1,\ldots,\bar{b}_{s}\in \bar{\mathfrak{B}}_{W(\Lambda)}$ 
such that vectors $\bar{b}_{1}v_{\Lambda}$, \ldots, $\bar{b}_{s}v_{\Lambda}$
span $W(\Lambda)_{\mu}$. Furthermore, the vectors $b_{1}v_{\Lambda}$, \ldots, $b_{s}v_{\Lambda}$,
where  $b_1,\ldots,b_{s}\in \mathfrak{B}_{W(\Lambda)}$, are linearly independent
elements of $W(\Lambda)_{\mu}$. We conclude that the vectors $\bar{b}_{1}v_{\Lambda}$, \ldots, $\bar{b}_{s}v_{\Lambda}$
are linearly independent and that the vectors $b_{1}v_{\Lambda}$, \ldots, $b_{s}v_{\Lambda}$ span $W(\Lambda)_{\mu}$.
The theorem now follows.
\end{dok}


\section{Quasi-particles of type 1 revisited}

For an arbitrary  vector space $V$  set 
$$\mathcal{E}(V):=\om(V,V((z))).$$
We recall two definitions from \cite{Li}.

\begin{defn}\label{quasi:comp}
An ordered sequence $(a_1(z),a_2(z),\ldots,a_{m}(z))$ in $\mathcal{E}(V)$ is said to be quasi compatible if there exist a nonzero polynomial 
$p(z_1,z_2)\in\mathbb{C}(q)[z_1,z_2]$
such that
$$\left(\prod_{r=1}^{m-1}\prod_{s=r+1}^{m}p(z_r,z_s)\right)a_1(z_1)a_2(z_2)\cdots a_{m}(z_m)\in\om (V,V((z_1,\ldots,z_m))).$$
\end{defn}

\begin{defn}
Let $(a(z),b(z))$ be a quasi compatible (ordered) pair in $\mathcal{E}(V)$. For $\alpha\in\mathbb{C}(q)\setminus\left\{0\right\}$, $l\in\mathbb{Z}$,
 we define
$a(z)_{(\alpha,l)}b(z)\in (\ndo V)[[z^{\pm 1}]]$ in terms of generating function
$$Y_{\mathcal{E}}^{(\alpha)}(a(z),z_0)b(z)=\sum_{l\in\mathbb{Z}} (a(z)_{(\alpha,l)}b(z))z_{0}^{-l-1}\in (\ndo V)[[z_{0}^{\pm 1}, z^{\pm 1}]]$$
by
\begin{align*}
& Y_{\mathcal{E}}^{(\alpha)}(a(z),z_0)b(z)
=\iota_{z,z_0}\left(p(z_0 +\alpha z,z)^{-1}\right)(p(z_1,z)a(z_1)b(z))\bigl. \bigr|_{z_1 =\alpha z+z_{0}},
\end{align*}
where $p(z_1,z_2)\in\mathbb{C}(q)[z_1,z_2]$, $p(z_1,z_2)\neq 0$, is any polynomial such that
\begin{equation}\label{poly7}
p(z_1,z_2)a(z_1)b(z_2)\in \om(V,V((z_1,z_2))).
\end{equation}
\end{defn}

In the rest of this section the parameter $\alpha$ will be equal to $1$ so we will omit it and write
$$Y_{\mathcal{E}}(a(z),z_0)b(z)=\sum_{l\in\mathbb{Z}} (a(z)_{l}b(z))z_{0}^{-l-1}.$$

The following lemma has a straightforward proof and, therefore, we skip it. 

\begin{lem}\label{l_poly}
On every integrable highest weight module $V$ we have 
\begin{align*}
&\left(\prod_{s=1}^{m+1}\left(z_1-q^{2(s-m)}z\right)\right)x_{\alpha_{i}}^{+}(z_1) x_{(m+1)\alpha_{i}}^{+}(q^{-2m}z)\\
&\hspace{40pt}=\left(\prod_{s=1}^{m+1}\left(q^{2}z_1-zq^{2(s-m-1)}\right)\right)x_{(m+1)\alpha_{i}}^{+}(q^{-2m}z)x_{\alpha_{i}}^{+}(z_1)\nonumber
\end{align*}
for $m\in\mathbb{Z}_{\geq 0}$ and $i=1,2,\ldots,n$.
In particular, we have 
$$\left(\prod_{s=1}^{m+1}\left(z_1-q^{2(s-m)}z\right)\right)x_{\alpha_{i}}^{+}(z_1) x_{(m+1)\alpha_{i}}^{+}(q^{-2m}z)\in\om(V,V((z_1,z))).$$
\end{lem}

The quasi compatiblity of an ordered pair $(x_{\alpha_{i}}^{+}(z), x_{(m+1)\alpha_{i}}^{+}(q^{-2m}z))$, established by the above lemma, will
be 
used in the proof of the following proposition.  

\begin{pro}
On every integrable highest weight module $V$ we have
\begin{align}\label{hsli}
&x_{\alpha_{i}}^{+}(zq^{-2m})_{-1}\Bigl(\ldots \left(x_{\alpha_{i}}^{+}(zq^{-4})_{-1}(x_{\alpha_{i}}^{+}(zq^{-2})_{-1}x_{\alpha_{i}}^{+}(z))\right)\ldots\Bigr)\\
&\hspace{30pt} =\left(\prod_{r=1}^{m}\prod_{s=r+1}^{m+1}\left(1-q^{2(s-r)+2}\right)\right)x_{(m+1)\alpha_{i}}^{+}(q^{-2m}z).\nonumber
\end{align}
for $i=1,2,\ldots,n$ and $m\in\mathbb{Z}_{>0}$.
\end{pro}

\begin{dok}
The proposition is proved by induction. The basis of induction is a formula
\begin{equation*}
x_{\alpha_{i}}^{+}(zq^{-2})_{-1}x_{\alpha_{i}}^{+}(z)=(1-q^{4})x_{2\alpha_{i}}^{+}(q^{-2}z),
\end{equation*}
that can be easily verified by a direct calculation. We assume that (\ref{hsli}) holds. Next, in the step of induction we use the polynomial
$$p(z_1,z)=\prod_{s=1}^{m+1}\left(z_1-q^{2(s-m)}z\right),$$
obtained in Lemma \ref{l_poly}, in order to calculate
$$x_{\alpha_{i}}^{+}(zq^{-2(m+1)})_{-1}x_{(m+1)\alpha_{i}}^{+}(q^{-2m}z),$$ 
thus proving the proposition.
\end{dok}

\begin{rem}
The left hand side of  equality (\ref{hsli}) is well defined
on every restricted module of $U_{q}(\widehat{\mathfrak{sl}}_{n+1})$ (see \cite{Li2})
so we can employ it in order to generalize our Definition (\ref{type_1_op.}) of the operator $x_{(m+1)\alpha_{i}}^{+}(z)$ to 
restricted modules.

Therefore,  on every restricted module we can define
\begin{align*}
x_{(m+1)\alpha_{i}}^{+}(z):=&\left(\prod_{r=1}^{m}\prod_{s=r+1}^{m+1}\frac{1}{1-q^{2(s-r)+2}}\right)\\
&\hspace{20pt}\cdot x_{\alpha_{i}}^{+}(z)_{-1}\Bigl(\ldots \left(x_{\alpha_{i}}^{+}(zq^{2(m-2)})_{-1}(x_{\alpha_{i}}^{+}(zq^{2(m-1)})_{-1}x_{\alpha_{i}}^{+}(zq^{2m}))\right)\ldots\Bigr)
\end{align*}
for $m\in\mathbb{Z}_{>0}$, thus generalizing (\ref{type_1_op.}) and Definition \ref{type_1_def}.
Naturally, for $m=0$  we set
$x_{1\alpha_{i}}^{+}(z):=x_{\alpha_{i}}^{+}(z)$.
\end{rem}


\section*{Acknowledgement}

Results of this paper are part of author's Ph.D. dissertation.
I would like to express my  gratitude to my Ph.D. advisor,
prof. M. Primc for  valuable guidance and   willingness to give his time so generously.


\end{document}